\documentclass[a4paper,12pt]{article}
\usepackage{chngcntr}
\usepackage{multirow,hhline}
\usepackage[margin=1.2in]{geometry}
\pdfoutput=1
\usepackage{mathtools}

\usepackage{amsfonts,amssymb, bm, amsthm}
\usepackage{graphicx,color,epstopdf}
\usepackage{showlabels} 
\allowdisplaybreaks
\newtheorem{lemma}{Lemma}
\newtheorem{theorem}{Theorem}
\newtheorem{remark}{Remark}

\def\XXint#1#2#3{{\setbox0=\hbox{$#1{#2#3}{\int}$}
    \vcenter{\hbox{$#2#3$}}\kern-.5\wd0}}

\def\hG {\hat{G}}

\def\tx {\tilde{x}}
\def\tk {\tilde{k}}

\def\tomg {\tilde{\omega}}

\def\hu {\hat{u}}
\def\hf {\hat{f}}

\def\pd {{\partial}}

\def\hG {\hat{G}}
\def\hu {\hat{u}}
\def\hf {\hat{f}}

\def\bi{{\bf i}}

\begin{document}
\title{On exact truncation of backward waves in elastrodynamics}
\author{Wangtao Lu$^1$} \footnotetext[1]{Corresponding author. School of Mathematical
  Sciences, Zhejiang University, Hangzhou 310027, China. Email:
  wangtaolu@zju.edu.cn. This author is partially supported by NSFC Grant
  12174310 and by NSF of Zhejiang Province for Distinguished Young Scholars
  (LR21A010001). }

\maketitle
\begin{abstract}
  For elastic wave scattering problems in unbounded anisotropic media, the
  existence of backward waves makes classic truncation techniques fail
  completely. This paper is concerned with an exact truncation technique for
  terminating backward elastic waves. We derive a closed form of elastrodynamic
  Green's tensor based on the method of Fourier transform and design two
  fundamental principles to ensure its physical correctness. We present a
  rigorous theory to completely classify the propagation behavior of Green's
  tensor, thus proving a conjecture posed by B\'ecache, Fauqueux and Joly ({\bf
    J. Comp. Phys., 188, 2003}) regarding a necessary and sufficient condition of
  the non-existence of backward waves. Using Green's tensor, we propose a new
  radiation condition to characterize anisotropic scattered waves at infinity.
  This leads to an exact transparent boundary condition (TBC) to truncate the
  unbounded domain, regardless the existence of backward waves or not. We
  develop a fast algorithm to evaluate Green's tensor and a high-accuracy scheme
  to discretize the TBC. A number of experiments are carried out to validate the
  correctness and efficiency of the new TBC.
\end{abstract}

\section{Introduction}
Acoustic, electromagnetic (optical), and elastic (seismic) waves form the
fundamental class of waves in real world, and have presented tremendous
applications in the past centuries, ranging from wireless communication, medical
imaging, to underwater and geological explorations, etc. Mathematically, the
propagation behavior of such three waves scattered by bounded obstacles are
governed by three well-known partial differential equations named after
Helmholtz, Maxwell and Navier, respectively. Specifically, if the propagation
medium is isotropic and homogeneous, all of the three scattering problems can be
formulated in close relation with Helmholtz's exterior problem,
\begin{align}
  \label{eq:helm}
  -\Delta u(x) - k^2 u(x) =& 0,\quad \mathbb{R}^d\backslash\overline{D},\\
  \label{eq:helm:bc}
  u =& f,\quad {\rm on}\quad\partial D,
\end{align}
where $d=2,3$ indicates the dimension of the problem, $x=(x_1,\cdots,x_d)$,
$\Delta = \sum_{j=1}^d\partial_{j}^2$ is the $d$-dimensional Laplacian with
$\partial_j=\frac{\partial}{\partial x_j}$, $u$ is a time-harmonic wave field,
$k>0$ is a constant depending on the medium property, $D$ represents a bounded
domain, and $f$ represents a given function.

This classic problem, firstly studied by Helmholtz in 1860, encountered an
essential difficulty: the solution may be non-unique. This wasn't resolved until,
in 1912, Sommerfeld in his pioneer work \cite{som12} firstly proved the
uniqueness by imposing an extra but necessary condition for $u$ at infinity:
\begin{equation}
  \label{eq:src}
  \lim_{r\to\infty}r^{(d-1)/2}\left( \frac{\partial u}{\partial r} - \bi k u\right) = 0,\quad r=|x|,\quad \bi = \sqrt{-1},
\end{equation}
uniformly w.r.t. all directions $x/r$. This condition, now
known as the famous Sommerfeld radiation condition (SRC), reveals an important
fact:
\begin{center}
  ``Any scattered wave must behave as an outgoing wave $r^{-(d-1)/2}e^{\bi k r}$ at infinity''.
\end{center}
It had subsequently been extended to the Silver-M\"{u}ller-Sommerfeld radiation
condition \cite{sil49, mul57} for Maxwell's equations and the
Kupradze-Sommerfeld radiation condition \cite{kup79} for Navier's equations. For
related well-posedness theories, we refer readers to \cite{colkre13, mcl00,
  kup79} and the references therein for details. Because of the significance of
the underlying motivation of SRC (\ref{eq:src}), we briefly review its
derivation first, according to the fabulous review paper \cite{sch92}.
Throughout the rest of the paper, we assume $d=2$ for simplicity.

\subsection{Motivation of SRC}
To derive (\ref{eq:src}), we require Green's function of (\ref{eq:helm}),
\begin{equation}
  \label{eq:phik}
  \Phi_k(x;y):=\frac{\bi}{4}H_0^{(1)}(k|x-y|) = \frac{\bi}{4\pi}\int_{-\infty}^{\infty}\frac{e^{\bi \sqrt{k^2-\xi^2}|x_2-y_2| + \bi \xi(x_1-y_1)}}{\sqrt{k^2-\xi^2}}d\xi,
\end{equation}
a.k.a. a fundamental solution, solving
\begin{equation}
  \label{eq:gov:Phi}
  -\Delta \Phi_k(x;y) - k^2 \Phi_k(x;y) = \delta(x-y),
\end{equation}
where $y=(y_1,y_2)$ denotes the source point and $H_m^{(1)}$ denotes the
first-kind Hankel function of $m$-th order.

Let $B(0,r)$
be a disk of sufficiently large radius $r$ and
$\Omega=B(0,r)\backslash\overline{D}$. Apply Green's third identity in $\Omega$,
\begin{equation}
  u(x) = \left[ \int_{\partial B(0,r)} + \int_{\partial D} \right]\left[\partial_{\nu(y)}\Phi_k(y;x) u(y) - \Phi_k(y;x)\partial_{\nu} u(y) \right] ds(y), x\in \Omega,
\end{equation}
where $\nu$ is the outer unit normal of $\partial\Omega$. Sommerfeld's original idea
is to propose certain condition to ensure the first integral over
$\partial B(0,r)$ vanishes as $r\to \infty$, i.e.,
\begin{equation}
  \label{eq:surc}
  \lim_{r\to\infty}\int_{\partial B(0,r)}\left[\partial_{\nu(y)}\Phi_k(y;x) u(y) - \Phi_k(y;x)\partial_{\nu} u(y) \right] ds(y) = 0,
\end{equation}
uniformly for $x$ in any bounded domain, since it directly yields Green's
exterior representation formula
\begin{equation}
  \label{eq:grep:helm}
  u(x) = \int_{\partial D}\left[\partial_{\nu(y)}\Phi_k(y;x) u(y) - \Phi_k(y;x)\partial_{\nu} u(y) \right] ds(y), x\in \mathbb{R}^2\backslash\overline{D}.
\end{equation}
By
\[
  \partial^j_{|y|}\Phi_k(x;y) = (\bi k)^j\sqrt{\frac{2}{\pi k |y|}}e^{\bi k
    |y|}+ {\cal O}(|y|^{-3/2}), j=0,1,2, \quad |y|\to\infty,
\]
uniformly for $x$ in any bounded domain,
(\ref{eq:surc}) becomes
\begin{equation}
  \lim_{r\to\infty}e^{\bi k r}\int_{0}^{2\pi}\sqrt{r}\left[\bi ku - \partial_{r} u \right]d\theta = 0,
\end{equation}
which immediately motivates SRC (\ref{eq:src}). Roughly speaking, SRC
(\ref{eq:src}) is an elegant way of stating the equivalent but unpleasant
condition (\ref{eq:surc}), which will be referred to as Sommerfeld's unpleasant
radiation condition (SURC).

\subsection{Truncation techniques based on SRC}
With SRC condition (\ref{eq:src}), the well-posedness of Helmholtz's exterior 
problem (\ref{eq:helm}) had been settled \cite{colkre13}. However, it remains challenging to
numerically compute $u$ due to the unbounded computational domain. Nevertheless,
SRC motivates three widely used truncation techniques.

Due to the special form of Laplacian in polar coordinates (or spherical
coordinates if $d=3$), the method of variable separation is applicable in $r\geq
r_0$, for some sufficiently large $r_0>0$ so that $D\subset B(0,r_0)$, to
express
\[
  u(r,\theta) = \sum_{n=0}^{\infty} \left[a_n H_n^{(1)}(k r)+ b_n H_n^{(2)}(k
    r)\right]e^{\bi m\theta},\quad r\geq r_0.
\]
SRC implies that $b_n=0$ for $n\in\mathbb{N}$ since $H_n^{(2)}$ is not outgoing
and does not satisfy (\ref{eq:src}). Then, the above induces a
variable-separation based transparent boundary condition (TBC)
\begin{equation}
  \label{eq:tbc1:helm}
  \partial_r u = \Lambda u,\quad {\rm on}\quad r = r_0,
\end{equation}
for a bounded operator $\Lambda: H^{1/2}(\partial B(0,r_0))\to H^{-1/2}(\partial
B(0,r_0))$. The second, more general approach is based on Green's representation
formula (\ref{eq:grep:helm}). Letting $x$ approaching $\partial D$ and using
standard jump relations yield a boundary-integral-equation (BIE) based TBC
\begin{equation}
  \label{eq:tbc2:helm}
  {\cal A} u + {\cal B}\partial_{\nu} u = 0,\quad{\rm on}\quad \partial D,
\end{equation}
for two bounded integral operators ${\cal A}: H^{1/2}(\partial D)\to
H^{1/2}(\partial D)$ and ${\cal B}: H^{-1/2}(\partial D) \to H^{1/2}(\partial
D)$. For related works using analogous TBC techniques, readers are referred to
\cite{ammbaowoo00, colkre13, gerdem96, mon03, kre14, baoxuyin17} and the references
therein.

The previous two approaches induce nonlocal boundary conditions to truncate the
unbounded domain. The third, more attractive approach is to develop a local
boundary condition by introducing an artificial absorbing medium, the so-called
perfectly matched layer (PML), originally developed by Berenger \cite{ber94} in
1994 for time-domain Maxwell's equations. Mathematically, PML introduces the
following complexified coordinates \cite{chewee94}
\begin{equation}
  \label{eq:pml:xi}
  \tx_i = x_i + \bi\int_{0}^{x_i} \sigma_i(t) dt, i=1,2,
\end{equation}
where $\sigma_i(t)=0$ for $|t|\leq L_i$, and are nonzero for $L_i\leq |t|\leq
L_i+d_i$. The regions of nonzero $\sigma_i$ constitute the PML region. Any
scattered wave $u$ satisfying SRC (\ref{eq:src}) and hence (\ref{eq:grep:helm}) shall be
complexified to
\[
  u(\tx) =\int_{\partial D}\left[\partial_{\nu(y)}\Phi_k(y;\tx) u(y) -
    \Phi_k(y;\tx)\partial_{\nu} u(y) \right] ds(y),\quad x\in \mathbb{R}^2\backslash\overline{D}. 
\]
As has been indicated in \cite{chewu03}, the outgoing behavior of $\Phi_k$
making its complexification $\Phi_k(y;\tx)$ decays exponentially in the PML and
so does $u(\tx)$. Therefore, it is spectrally accurate to impose a zero
Dirichlet boundary condition on the PML boundary, i.e.,
\begin{equation}
  \label{eq:pmlbc}
  u(\tx) = 0,\quad {\rm on}\quad \{x:|x_1| = L_1+d_1{\rm\ or\ } |x_2|=L_2+d_2\}.
\end{equation}
Numerically computing $u(\tx)$ recovers $u(x)$ with high accuracy outside the PML for $\tx=x$. Due
to its simplicity and nearly zero reflection for outgoing waves, PML has, since
its development, been widely used in the simulation of a variety of wave
propagation problems \cite{tafhug05, rodged00}. Readers are referred to
\cite{chewu03,chexiazha16,chezhe10,chezhe17, lassom01, lulaiwu19} for the
spectral convergence theories of the PML truncation for various scattering
problems, and to \cite{engmaj77} for an approximate but effective absorbing
boundary condition due to Engquist and Majda. Besides, the author and his
collaborators have recently developed a hybrid truncation technique that
combines PML and a BIE-based TBC technique for scattering problems in
complicated structures \cite{luluqia18,lu21,gaolu22,yuhulurat22}.

\subsection{Elastrodynamics: failure of SRC}
As has been seen, the related truncation techniques rely mostly on the outgoing
behavior of the scattered field due to the isotropicity of the media.
Unfortunately, such a significant prerequisite breaks down in anisotropic media,
due to the possible existence of backward waves. Physically speaking, a backward
wave refers to a wave with an always outgoing group velocity, the physically
correct direction of energy transport \cite{joajohwinmea08}, propagates incoming
from the infinity. For Helmholtz's scattering problem (\ref{eq:helm}), there is
no backward wave since a plane wave $e^{\bi (v_1 x_1 + v_2 x_2)}$ of directional
vector $v=(v_1,v_2)$ satisfies $k= \sqrt{v_1^2 + v_2^2}$ so that the group
velocity $\nabla_vk$ always coincides with the propagation direction $v$. This
is unfortunately not true in elastrodynamics with anisotropic media commonly
encountered in practice \cite{car15}. Consequently, it becomes a fundamental
task to precisely characterize and accurately simulate the propagation
behavior of elastic waves in such media. This is the main objective of the
current paper.

Let us start with a Dirichlet-type exterior problem. Now, $D$ denotes a bounded
rigid body with its complement $\mathbb{R}^2\backslash\overline{D}$ filled by a
homogeneous but anisotropic elastic medium. A time-harmonic scattered wave, of
time dependence $e^{-\bi\omega t}$,
$u=[u_1,u_2]^{T}$ due to an incident wave $u^{\rm inc}=[u^{\rm inc}_1, u^{\rm
  inc}_2]^{T}$ imping upon $D$ is then governed by Navier's equation
\begin{align}
  \label{eq:gov:2d}
  -\nabla\cdot\sigma(u) - \rho\omega^2u =& 0,\quad\mathbb{R}^2\backslash\overline{D},\\
  \label{eq:dir:bc}
  u(x) =& -u^{\rm inc},\quad{\rm on}\quad\partial D,
\end{align}
where (\ref{eq:dir:bc}) is due to the zero total wave $u^{\rm tot}=u + u^{\rm
  inc}$ on $\partial D$, $\nabla = [\pd_1,\pd_2]^{T}$, $\omega$ is the frequency of the incident
wave, $\sigma$ is the stress tensor satisfying
\begin{equation}
  \sigma = C\epsilon(u),
\end{equation}
$C=[C_{ijkl}]_{i,j,k,l=1}^2$ is a fourth-order tensor,
the strain tensor $\epsilon(u)=[\epsilon_{ij}(u)]$ with
\[
  \epsilon_{ij}(u) = \frac{1}{2}(\pd_i u_j+\pd_j u_i),
\]
and
\begin{equation}
  \nabla\cdot\sigma = \left[
    \begin{array}{l}
      \pd_1\sigma_{11}+\pd_2\sigma_{12}\\ 
      \pd_1\sigma_{12}+\pd_2\sigma_{22}
      \end{array}
  \right].
\end{equation}
Using the Voigt notation \cite{car15}, $C$ can be represented by a $3\times 3$
positive symmetric matrix denoted by 
\[
  C = \left[
    \begin{array}{ccc}
      c_{11} &c_{12}& c_{13}\\
      c_{12} &c_{22}& c_{23}\\
      c_{13} & c_{23} & c_{33}
    \end{array}
  \right],
\]
and (\ref{eq:gov:2d}) becomes
\begin{equation}
  \label{eq:gov}
  -\pd_1(A_{11}\pd_1 u + A_{12}\pd_2 u) -\pd_2(A_{21}\pd_1 u + A_{22}\pd_2 u) -\rho\omega^2 u = 0,\quad {\rm in}\quad \mathbb{R}^2\backslash\overline{D},
\end{equation}
where
\[
  A_{11} = \left[
    \begin{matrix}
      c_{11} & c_{13}\\
      c_{13} & c_{33}
    \end{matrix}
  \right], \quad  A_{12} = A_{21}^{T}=\left[
    \begin{matrix}
     c_{13} & c_{12} \\
     c_{33} & c_{23}
    \end{matrix}
  \right],\quad A_{22} = \left[
    \begin{matrix}
      c_{33} & c_{23}\\
      c_{23} & c_{22}
    \end{matrix}
  \right]. 
\]
To illustrate the basic idea, we consider orthotropic media which assume
$c_{13}=c_{23}=0$ throughout this paper; the case when $c_{13}\neq 0$ or
$c_{23}\neq 0$ can be analyzed with no technical difficulty. Thus, the nonzero
elements $c_{ij}$ must satisfy
\[
(C_0)\ \left\{
  \begin{array}{ll}
    (i) &c_{11},c_{22},c_{33}>0,c_{11}c_{22}-c_{12}^2>0;\\
    (ii) &c_{12}+c_{33}\neq 0,
  \end{array}
  \right.  
\]
where $(C_0)(ii)$ is to avoid the decoupling of $u_1$ and $u_2$.

Similar to Helmholtz's problem (\ref{eq:helm}), Navier's problem (\ref{eq:gov})
necessarily requires an accurate radiation condition at infinity, which remains
absent so far to the author's best knowledge. Even worse, the failure of SRC
makes the extension of the aforementioned truncation techniques precarious.
Specifically, the first technique via a variable-separation based TBC like
(\ref{eq:tbc1:helm}) fails as variable separations may not work in the first
place, let alone the invalidity of SRC. The second technique based on a
BIE-based TBC like (\ref{eq:tbc2:helm}) now faces a new challenging: Green's
tensor. A closed-form of Green's tensor based on $d$-dimensional Fourier
transform is trivial \cite{liulam96}. However, it remains unclear whether it
reveals the true radiation behavior of the scattered field, as Green's tensor
needs not to be unique.

The third technique based on PML now becomes more subtle. B\'{e}cache {\it et
  al.} in \cite{becbonfliton22} rigorously justified the instability of PML in
elastrodynamics due to the existence of backward waves, which blow up in the
PML. They further derived a sufficient condition to ensure the stability of PML,
and conjectured its necessity, the proof of which is incomplete yet. To try to
resolve the instability, Appel\"{o} and Kreiss \cite{appkre06} and Savadatti and
Guddati \cite{savgud12a,savgud12b} subsequently designed new types of PMLs,
which applied only for a limited class of elastic media. Recently, B\'{e}cache\
{\it et al.} \cite{becbonfliton22} developed a half-space matching method to
solve Navier's problem (\ref{eq:gov}). They applied the simple one-dimensional
Fourier transform along four straight lines surrounding the obstacle $D$,
proposed rules to capture backward waves, and successfully developed a
truncation technique by solving four half-space problems.

\subsection{Methodology: opportunity of SURC}
The methodology of the current work closely follows Sommerfeld's philosophy for
Helmholtz's problem (\ref{eq:helm}). We believe that SURC (\ref{eq:surc}) can be
generalized to prescribe the asymptotic behavior of acoustic, electromagnetic,
and elastic scattered waves in homogeneous but whatever anisotropic media, if
the associated Green tensor/function is correctly selected! Meanwhile, the
selecting principle should obey the same rule in choosing $\Phi_k$ for
Helmholtz's problem (\ref{eq:helm}), as discussed below.

The integral form of $\Phi_k(x;y)$ in (\ref{eq:phik}), a.k.a. Sommerfeld's
integral, is derived by Fourier transforming (\ref{eq:gov:Phi}) w.r.t. the
single variable $x_1$ (generally, the $d-1$ variables $x_1,\cdots, x_{d-1}$),
matching the Fourier transform of $\Phi_k$ for the other variable $x_2$ at
$x_2=y_2$, and an inverse Fourier transform finally. Then, the radiation
behavior of $\Phi_k$ is completely determined by a proper choice of branch cut
of $\sqrt{k^2-\xi^2}$; choosing a wrong branch cut could lead to incoming
Green's function $\overline{\Phi_k}(x;y)$. This is completely resolved by
replacing the real path with Sommerfeld's integral path (SIP), a complex path
going from an infinity in Quadrant II to the origin, and then to another
infinity in Quadrant IV. Now $\sqrt{k^2-\xi^2}$ is always complex (except at
$\xi=0$) and the basic rule is simply letting ${\rm Im}(\sqrt{k^2-\xi^2})\geq 0$
to ensure the boundedness of the integrand in (\ref{eq:phik}). The above
constitutes the basic idea of selecting proper Green's tensor for Navier's
problem (\ref{eq:gov}). Upon an analog of SURC (\ref{eq:surc}) and a carefully
selected Green tensor, an analog of Green's exterior representation formula
(\ref{eq:grep:helm}) for Navier's problem can be derived. Then, the second TBC
truncation technique becomes straightforwardly extendable, as it does not care at
all about the outgoing or incoming behavior of $u$!

The rest of this paper is organized as follows. In section 2, we derive Green's
tensor for Navier's problem (\ref{eq:gov}). In section 3, we present a rigorous
theory to completely classify the propagation behavior of Green's tensor for any
tensor $C$ satisfying $(C_0)$, thus providing an affirmative answer to the
aforementioned conjecture in \cite{becfaujol03}. In section 4, a fast algorithm
is proposed to evaluate Green's tensor and its derivatives. In section 5, we
propose an analog of SURC (\ref{eq:surc}) for Navier's problem (\ref{eq:gov}),
leading to a new TBC analogous to (\ref{eq:tbc2:helm}), and carry out a number
of numerical experiments to justify the correctness of the TBC. In section 6, we
study two typical scattering problems and demonstrate the efficiency of the new
TBC by a number of examples. In section 7, we conclude this paper by raising
some significant problems.

\section{Elastrodynamic Green's tensor}
The two-dimensional elastrodynamic Green tensor
\[
  G(x;y) = \left[
    \begin{matrix}
     G_{11}(x;y) & G_{12}(x;y)\\ 
     G_{21}(x;y) & G_{22}(x;y)\\ 
    \end{matrix}
  \right]
\]
excited by a source point $y\in\mathbb{R}^2$, is defined as a solution of
\begin{equation}
  \label{eq:gov:G}
  -\left[\pd_1(A_{11}\pd_1  + A_{12}\pd_2 ) + \pd_2(A_{21}\pd_1  + A_{22}\pd_2 ) + \rho\omega^2  \right] G(x;y) = I\delta(x-y),
\end{equation}
where $I$ denotes the $2\times 2$ identity matrix. This problem has
infinitely many solutions. In this section, we propose two fundamental principles to
single out a unique and physically reasonable solution $G(x;y)$, and derive its
closed form for any tensor $C$ satisfying $(C_0)$ by the method of Fourier
transform.

For simplicity, we assume $y=0$ for the moment, and write $G(x)$ short for
$G(x;0)$. Consider $G_1=[G_{11}, G_{21}]^{T}$ first. For a generic function
$f(x)$, let
\[
  \hf(x_2;\xi) = \int_{-\infty}^{\infty} f(x) e^{-\bi \xi x_1} dx_1
\]
be the $x_1$-Fourier transform of $f$. Then, $x_1$-Fourier transforming the first
column of (\ref{eq:gov:G}) yields
\begin{align}
  \label{eq:gov:1:fre}
  \bi\xi (c_{12}+c_{33})\frac{d}{dx_2}\hG_{21} + c_{33}\frac{d^2}{dx^2_2}\hG_{11} +(\rho\omega^2-c_{11}\xi^2)\hG_{11}=&-\delta(x_2),\\
  \label{eq:gov:2:fre}
  \bi\xi(c_{12}+c_{33})\frac{d}{dx_2}\hG_{11} + c_{22}\frac{d^2}{dx^2_2}\hG_{21} +(\rho\omega^2-c_{33}\xi^2)\hG_{21}=&0.
\end{align}
This implies the following continuity conditions
\begin{align}
  \label{eq:cond:1}
  [\hG_{11}]_{x_2=0} =& 0,\quad [\hG'_{11}]_{x_2=0} = -c_{33}^{-1},\\
  \label{eq:cond:2}
  [\hG_{21}]_{x_2=0} =& 0,\quad [\hG'_{21}]_{x_2=0} = 0,
\end{align}
where $[\cdot]_{x_2=0}$ indicates the jump of the quantity from $x_2=0+$ to $x_2=0-$.

For $x_2\neq 0$, we obtain
\begin{align}
  \label{eq:gov:12}
\left[c_{33}\frac{d^2}{dx^2_2}+(\rho\omega^2-c_{11}\xi^2)  \right]\hG_{11} =-\bi\xi (c_{12}+c_{33})\frac{d}{dx_2}\hG_{21},\\ 
  \label{eq:gov:21}
\left[c_{22}\frac{d^2}{dx^2_2}+(\rho\omega^2-c_{33}\xi^2)  \right]\hG_{21} =-\bi\xi (c_{12}+c_{33})\frac{d}{dx_2}\hu_{11}. 
\end{align}
Thus, $\hG_{j1}$ is governed by the following characteristic equation
\begin{align}
  \label{eq:cha}
  c_{33}c_{22}\hG_{j1}^{(4)} + \left[  \rho\omega^2\beta_0 +\xi^2\beta_1\right]\hG_{j1}^{(2)}+ (\rho\omega^2-c_{11}\xi^2)(\rho\omega^2-c_{33}\xi^2) \hG_{j1}=0,
\end{align}
for $j=1,2$, where
\begin{align}
  \label{eq:b0}
  \beta_0=&c_{22}+c_{33},\\
  \label{eq:b1}
  \beta_1=&c_{12}^2+2c_{12}c_{33}-c_{11}c_{22}.
\end{align}
A simple derivation shows that 
\[
  \hG_{j1}(x_2;\xi) = \left\{
    \begin{array}{lc}
    p^+_{j}(\xi)e^{\bi\mu_+(\xi^2) x_2} + q^+_{j}(\xi)e^{\bi\mu_-(\xi^2) x_2},\quad x_2>0;\\
    p^-_{j}(\xi)e^{-\bi\mu_+(\xi^2) x_2} + q^-_{j}(\xi)e^{-\bi\mu_-(\xi^2) x_2},\quad x_2<0,
    \end{array}
  \right.
\]
for $j=1,2$, where $p_j^\pm$ and $q_j^\pm$ are unknowns to be determined and
$\mu_\pm^2$ are roots of 
\begin{equation}
  \label{eq:cha:mu}
c_{33}c_{22}\left( \frac{\eta}{\rho\omega^2} \right)^2 - \left(  \beta_0 +\frac{\xi^2}{\rho\omega^2}\beta_1\right)\frac{\eta}{\rho\omega^2} + \left(1-c_{11}\frac{\xi^2}{\rho\omega^2}\right)\left(1-c_{33}\frac{\xi^2}{\rho\omega^2}\right)=0.
\end{equation}
Thus,
\begin{align}
  \label{eq:mupm}
  \mu_\pm^2(\xi^2) =& \frac{\left[  \rho\omega^2\beta_0 +\xi^2\beta_1\right]\pm\sqrt{\Delta(\xi^2)} }{2c_{33}c_{22}},\\
  \label{eq:Delta}
  \Delta(\xi^2) 
  =&\rho^2\omega^4\alpha_0  + 2\rho\omega^2\xi^2\alpha_1 + \xi^4\alpha_2,
\end{align}
where $\alpha_0 =(c_{22}-c_{33})^2$ and
\begin{align}
  \alpha_1
  =& (c_{22}+c_{33})(c_{12}+c_{33})^2 - (c_{22}-c_{33})(c_{11}c_{22}-c_{33}^2),\\
  \label{eq:alpha2}
  \alpha_2
  =&(c_{11}c_{22}-c_{12}^2)(c_{11}c_{22}-(c_{12}+2c_{33})^2). 
\end{align}

In the above, three special square-root functions $\sqrt{\Delta}$ and
$\sqrt{\mu_\pm^2}$ appear. To correctly capture the propagation behavior of $G$,
we propose two fundamental principles to select proper branch cuts for them. As
$c_{11}^{-1}$ and $c_{33}^{-1}$ are zeros of $\mu^2_\pm$ and hence branch points
of $\mu_\pm$, indicating that $\mu_\pm(\xi^2)$ may not be analytic for
$\xi\in\mathbb{R}$. To resolve this, we make use of the aforementioned SIP. Let
\begin{align}
  \chi(t) = \left\{
  \begin{array}{lc}
    t, &0\leq t\leq 1;\\
    1,&{\rm otherwise}.
  \end{array}
  \right.
\end{align}
The following path
\begin{equation}
  \label{eq:z:sip}
  \{z(\xi;\epsilon_0)= {\rm sgn}(\xi)\sqrt{\xi^2-\epsilon\bi}: \epsilon(\xi^2)=\epsilon^2_0\chi(\xi^2/\epsilon_0), \xi\in\mathbb{R},0<\epsilon_0\ll 1\}
\end{equation}
in the complex plane is a typical SIP and shall be used throughout this paper
unless otherwise indicated (c.f. Remark~\ref{rem:gensip}),
where ${\rm sgn}$ is the standard sign function and, for simplicity, the
arguments of $z$ and $\epsilon$ shall be suppressed frequently. Unless otherwise indicated, the
branch cut of a usual square-root function shall be chosen as the negative real
axis so that its argument is limited to $(-\frac{\pi}{2}, \frac{\pi}{2}]$. Thus,
SIP goes from $-\infty+0\bi$ to the origin in Quadrant II, and then to
$+\infty-0\bi$ in Quadrant IV; here $\pm 0\bi$ indicates the path approaches the
real axis at infinity from above and below, respectively. The first principle is
\begin{center}
  (P1)\quad\quad$\mu_\pm$ must be analytic along SIP for $0<\epsilon_0\ll 1$;
\end{center}
this is based on that $G$ is analytic for $x_2\neq 0$. Based on
(P1), we define
\begin{equation}
  \label{eq:mu:def}
  \mu_\pm(\xi^2) = \lim_{\epsilon_0\to 0^+}\mu_\pm(z^2(\xi;\epsilon_0)).
\end{equation}
The second principle is 
\begin{center}
  (P2)\quad\quad${\rm Im}(\mu_\pm(\xi^2))\geq 0,\quad\forall\xi\in\mathbb{R}$;
\end{center}
otherwise, $\hG_{1j}$ blows up at $x_2=\pm\infty$. In section 3, we shall
determine the branch cuts of the three functions to realize (P1) and (P2). We
first make some preliminary remarks regarding the choice of branch cuts for the
three functions.
\begin{remark}
The branch cut of $\sqrt{\Delta}$ is not so essential as it just affects
the definitions of $\mu_\pm^2$ in (\ref{eq:mupm}); we only require that
$\sqrt{\Delta}$ is analytic along the SIP to protect (P1). To achieve (P2), we
may directly adopt the positive real axis as the branch cuts of
$\sqrt{\mu_\pm^2}$ so that $\arg\mu_\pm\in[0,\pi)$. However, this cannot
guarantee (P1). 
\end{remark}
From the continuity conditions (\ref{eq:cond:1}) and (\ref{eq:cond:2}) and the
two governing equations (\ref{eq:gov:12}) and (\ref{eq:gov:21}), the eight
unknowns $p_j^\pm$ and $q_j^\pm$ can be derived on the SIP to bypass zeros of
$\Delta$ and $\mu_\pm^2$. They are given by
\begin{align}
  \label{eq:p2}
  p_2^+=-p_2^-=-q_2^+=q_2^- = \frac{-\bi z(c_{12}+c_{33})c_{33}^{-1}}{2c_{22}(\mu_+^2-\mu_-^2)} = \frac{-\bi z(c_{12}+c_{33})}{2\sqrt{\Delta(z^2)}},
\end{align}
and
\begin{align}
  \label{eq:p1}
  p_1^+=&p_1^- =-\bi\frac{\left[-c_{22}\mu_+^2(z^2)+(\rho\omega^2-c_{33}z^2)  \right]}{2\sqrt{\Delta(z^2)}\mu_+(z^2)},\\
  \label{eq:q1}
  q_1^+=&q_1^- =\bi\frac{\left[-c_{22}\mu_-^2(z^2)+(\rho\omega^2-c_{33}z^2)  \right]}{2\sqrt{\Delta(z^2)}\mu_-(z^2)}.
\end{align}
Then, $G_1$ can be obtained via the following inverse Fourier transform,
\[
  G_1(x) = \frac{1}{2\pi}\int_{\mathbb{R}} \hG_1(x;\xi) e^{\bi\xi x_1}d\xi=
  \frac{1}{2\pi}\int_{\rm SIP} \hG_1(x;z) e^{\bi z x_1}dz,
\]
where the second equality holds by Cauchy's theorem and by the two principles
(P1) and (P2). Consequently, by a simple translation, one obtains
$G_1(x;y)=[G_{11}(x;y);G_{21}(x;y)]^{T}$ given by
\begin{align}
  \label{eq:green11}
  G_{11}(x;y) =&\frac{1}{2\pi}\int_{\mathbb{R}}\left[ p^{x_1}_{1}(\xi)e^{\bi\mu_+(\xi^2) |x_2-y_2|} + q^{x_1}_{1}(\xi)e^{\bi\mu_-(\xi^2) |x_2-y_2|} \right]e^{\bi \xi |x_1-y_1|} d\xi,\\
  \label{eq:green21}
  G_{21}(x;y) =&\frac{-{\rm sgn}[(x_1-y_1)(x_2-y_2)]}{2\pi}\int_{\mathbb{R}}p^{x_1}_{2}(\xi)\left[e^{\bi\mu_+(\xi^2) |x_2-y_2|} - e^{\bi\mu_-(\xi^2) |x_2-y_2|}  \right]e^{\bi \xi |x_1-y_1|} d\xi,
\end{align}
where $p_1^{x_1}=p_1^+$, $p_2^{x_1}=p_2^+$, and $q_1^{x_1}=q_1^+$. Similarly, the second column of $G(x;y)$,
$G_2(x;y)=[G_{12}(x;y);G_{22}(x;y)]^{T}$ is derived to be
\begin{align}
  \label{eq:green12}
  G_{12}(x;y) =&-\frac{{\rm sgn}[(x_1-y_1)(x_2-y_2)]}{2\pi}\int_{\mathbb{R}}p^{x_2}_{1}(\xi)\left[ e^{\bi\mu_+(\xi^2) |x_2-y_2|} - e^{\bi\mu_-(\xi^2) |x_2-y_2|} \right]e^{\bi \xi |x_1-y_1|} d\xi,\\
  \label{eq:green22}
  G_{22}(x;y) =&\frac{1}{2\pi}\int_{\mathbb{R}}\left[p^{x_2}_{2}(\xi)e^{\bi\mu_+(\xi^2) |x_2-y_2|} + q^{x_2}_2(\xi)e^{\bi\mu_-(\xi^2) |x_2-y_2|}  \right]e^{\bi \xi |x_1-y_1|} d\xi,
\end{align}
where $p_1^{x_2}=p_2^{x_1}$, and
\begin{align}
  p_2^{x_2} =&-\bi\frac{-c_{33}\mu_+^2+(\rho\omega^2-c_{11}\xi^2)}{2\sqrt{\Delta(\xi^2)}\mu_+(\xi^2)},\\ 
  q_2^{x_2}=&\bi\frac{-c_{33}\mu_-^2+(\rho\omega^2-c_{11}\xi^2)}{2\sqrt{\Delta(\xi^2)}\mu_-(\xi^2)}.
\end{align}
Clearly, Green's tensor $G$ is symmetric, satisfies $G(x;y)=G(x-y)$, and obeys
the usual reciprocal relation
\[
  G(x;y) = G(y;x).
\]
To simplify the presentation, we shall assume $\rho\omega^2 = 1$ unless
otherwise indicated. For numerical purposes, we split $G$ into two parts as
\begin{equation}
  \label{eq:G:sp}
  G(x) = G^+(x) + G^-(x),
\end{equation}
with
\begin{align*}
  G^+(x) =\left[G^+_{ij}\right]=& \frac{1}{2\pi}\int_{\mathbb{R}}\left[
  \begin{matrix}
    p_1^{x_1}(\xi) & -{\rm sgn}(x_1x_2)p_2^{x_1}(\xi)\\
-{\rm sgn}(x_1x_2)p_2^{x_1}(\xi) & p_2^{x_2}(\xi)
  \end{matrix}
  \right]e^{\bi \mu_+|x_2| + \bi \xi |x_1|} d\xi,\\
  G^-(x) =\left[G^-_{ij}\right]=& \frac{1}{2\pi}\int_{\mathbb{R}}\left[
  \begin{matrix}
    q_1^{x_1}(\xi) & {\rm sgn}(x_1x_2)p_2^{x_1}(\xi)\\
{\rm sgn}(x_1x_2)p_1^{x_2}(\xi) & q_2^{x_2}(\xi)
  \end{matrix}
  \right]e^{\bi \mu_-|x_2| + \bi \xi |x_1|} d\xi.
\end{align*}

From the above, it can be seen that $G(x)$ consists of plane waves
\[
  P_\pm(x;\xi):=e^{\bi \mu_\pm |x_2| + \bi\xi x_1}, \xi\in\mathbb{R}.
\]
Clearly, (P2) ensures ${\rm Im}(\mu_\pm)\geq 0$ so that $P_\pm(x;\xi)$ is
bounded and will never blow up at $x=\infty$ for all $\xi\in\mathbb{R}$.
However, (P2) cannot determine the sign of ${\rm Re}(\mu_\pm)$, which determines
the propagation behavior of $P_\pm(x;\xi)$. Basically, there are two possible
situations:
\[
  P_\pm(x;\xi)\ {\rm is}
    \begin{cases}
      \textrm{outgoing along $x_2-$direction if} & {\rm Re}(\mu_\pm)\geq 0;\\ 
      \textrm{incoming along $x_2-$direction if} & {\rm Re}(\mu_\pm)\leq 0. 
    \end{cases}
\]
The outgoing/incoming waves $P_\pm(x;\xi)$ are referred to as the so-called
forward/backward propagating waves. Note that the ambiguity for ${\rm
  Re}(\mu_\pm)=0$ is safe here since it corresponds to the case when
$P_\pm(x;\xi)$ propagates purely along $x_1$-direction, and, if ${\rm
  Im}(\mu_\pm)>0$, decays along $x_2$-direction. In section 3, we shall
present a complete classification of Green's tensor $G$ for any tensor $C$
satisfying $(C_0)$ in the sense that the branch cuts of $\sqrt{\Delta}$ and
$\sqrt{\mu_\pm^2}$ as well as the propagation behavior of $P(x;\xi)$ are
completely revealed. Before that, we present a lemma that reveals the locations where $\mu_\pm^2(z^2)$ cross the real axis as $z$ travels along the SIP. They play a central role in
identifying the branch cuts of $\sqrt{\mu_\pm^2}$.
\begin{lemma}
  \label{lem:realcros}
  For $0<\epsilon_0\ll 1$, the two paths $\{\mu_\pm^2(z^2):\xi\in\mathbb{R}\}$
  in $\mathbb{C}$ together cross the real axis at most twice (Here and
  hereafter, we only count the number of nonzero values of $\xi^2$).
  Specifically,
\begin{itemize}
\item[(i).] If $\beta_1=0$, two crossing points occur at the same point
  \begin{equation}
    \label{eq:xi20}
    \xi_0^2 = \frac{c_{11}+c_{33}}{2c_{11}c_{33}}>0,
  \end{equation}
  with $\mu_+^2(\xi_0^2-\bi\epsilon)\mu_-^2(\xi_0^2-\bi\epsilon)<0$;
\end{itemize}
\item[(ii).] If $\beta_1\neq 0$, the crossing points, if any, are
\begin{equation}
  \label{eq:rootformula}
  \eta_0 = \frac{2c_{11}c_{33}\xi^2-(c_{11}+c_{33})}{\beta_1},\quad 
\end{equation}
with $\xi^2>0$ solving
\begin{equation}
  \label{eq:realroot}
  \alpha_2 \xi^4 + 2\alpha_1 \xi^2 -\frac{\gamma_2}{c_{11}c_{33}} + \beta_1^2\epsilon^2= 0,
\end{equation}
where
\begin{equation}
  \label{eq:g2}
  \gamma_2=\{(c_{12}+c_{33})^2 -
c_{11}(c_{22}-c_{33})\}\times\{(c_{12}+c_{33})^2 + c_{33}(c_{22}-c_{33})\}.
\end{equation}
\begin{proof}
  With $z^2=\xi^2-\epsilon\bi$ in place of $\xi^2$ in (\ref{eq:cha:mu}), the
  crossing points correspond to real roots of the following equation
\[
  \eta^2 -\frac{\beta_0 + (\xi^2-\epsilon\bi)\beta_1}{c_{33}c_{22}}\eta +
  \frac{(1- c_{11}(\xi^2-\epsilon\bi))(1-c_{33}(\xi^2-\epsilon\bi))}{c_{33}c_{22}} = 0.
\] 
Suppose it has a real root $\eta_0$. Equating the imaginary and real parts of the two
sides yields
\[
  \frac{\beta_1\epsilon}{c_{33}c_{22}}\eta_0 + \frac{(c_{11}+c_{33})\epsilon -
    2c_{11}c_{33}\epsilon \xi^2}{c_{33}c_{22}} = 0,
\]
and
\[
 \eta_0^2 -\frac{\beta_0 + \beta_1 \xi^2}{c_{33}c_{22}}\eta_0 + \frac{(1- c_{11}\xi^2)(1-c_{33}\xi^2)- c_{11}c_{33}\epsilon^2}{c_{33}c_{22}} = 0.
\]
If $\beta_1=0$, $\xi^2 = \xi_0^2$ so that $(1-c_{11}\xi_0^2)(1-c_{33}\xi_0^2)<0$
and $\mu_+^2(\xi_0^2-\bi\epsilon)\mu_-^2(\xi_0^2-\bi\epsilon)<0$. Suppose now
$\beta_1\neq 0$. For $\xi^2>0$ so that $\epsilon>0$, (\ref{eq:rootformula})
holds and hence $\xi^2>0$ must be a root of (\ref{eq:realroot}).
\end{proof}
\end{lemma}

\section{Classification of Green's tensor}
In this section, we determine the branch cuts of $\sqrt{\Delta}$ and
$\sqrt{\mu_\pm^2}$, and, as a by-product, give a complete classification of $G$
at infinity. We discuss $G$ at $x_2=\pm \infty$ only; for $G$ at
$x_1=\pm\infty$, one simply permutes $c_{11}$ and $c_{22}$ in $C$ to swap $x_1$
and $x_2$. It turns out that the sign of $\min_{\xi\in\mathbb{R}}\Delta(\xi^2)$
significantly affects the behavior of $G$ at infinity. Thus, we distinguish two
cases: $\min_{\xi\in\mathbb{R}}\Delta(\xi^2)\geq 0$ and
$\min_{\xi\in\mathbb{R}}\Delta(\xi^2)<0$.
\subsection{Purely outgoing Green's tensor for $\min_{\xi\in\mathbb{R}}\Delta(\xi^2)\geq 0$}
To ensure $\min\Delta(\xi^2)\geq 0$, B\'{e}cache, Fauqueux
and Joly \cite{becfaujol03} have found the following sufficient and necessary
condition: the tensor $C$ must satisfy
\begin{align*}
  {\rm (BFJ)}_{x_2}^0: \left\{  
  \begin{array}{l}
  (C_1)'\ \gamma_1:=\{(c_{12}+c_{33})^2 - c_{22}(c_{11}-c_{33})\}\times\{(c_{12}+c_{33})^2 + c_{33}(c_{11}-c_{33})\}\leq 0,\\
  (C_2)\ \left\{
  \begin{array}{ll}
    (i) &(c_{12}+2c_{33})^2\leq c_{11}c_{22}\\
    (ii) &(c_{12}+c_{33})^2<c_{11}c_{22} + c_{33}^2 
  \end{array}
  \right.,\\
    \textrm{one of} \left\{ 
    \begin{array}{l}
  (C_3)_1\quad (c_{12}+c_{33})^2\leq (c_{11}-c_{33})(c_{22}-c_{33}),\\
      (C_3)_2\quad (c_{22}+c_{33})(c_{12}+c_{33})^2\geq (c_{22}-c_{33})(c_{11}c_{22}-c_{33}^2).
    \end{array}
    \right.
  \end{array}
\right.
\end{align*}
\begin{remark}
  \label{rem:bfjx20}
  As indicated in \cite{becfaujol03}, $(C_1)'$ and $(C_2)(ii)$ are superfluous
  based on the observation that
  \begin{align*}
    \beta_1 =& c_{12}(c_{12}+2c_{33}) - c_{11}c_{22} \leq (|c_{12}|-\sqrt{c_{11}c_{22}})\sqrt{c_{11}c_{22}} < 0,\\ 
               0\leq& \Delta(-\beta_0\beta^{-1}_1) = -4\gamma_1c_{22}c_{33}\beta_1^{-2}.
  \end{align*}
  Nevertheless, we keep them here for completeness.
\end{remark}
We first figure out the branch cut of $\sqrt{\Delta(z^2)}$ for $z\in$SIP with
$\epsilon_0\ll 1$ under (BFJ)$_{x_2}^0$. By
\begin{equation}
  \label{eq:del:z2}
  \Delta(z^2) = \Delta(\xi^2) - \alpha_2\epsilon^2 - 2(\alpha_2\xi^2 + \alpha_1)\epsilon\bi.
\end{equation}
and by $\min\Delta(\xi^2)\geq 0$, it is straightforward to see that
$\Delta(z^2)$ never crosses the negative axis. Thus, the negative real axis can
serve as the branch cut of $\sqrt{\Delta}$ when $\min\Delta(\xi^2)>0$. Note that
if $\min\Delta(\xi^2)=0$, $\sqrt{\Delta}$ becomes entire.

\subsubsection{A strong condition}
We start from a strong condition that guarantees no backward waves.
It states as follows:

\vspace{0.3cm}

\noindent(OC)$_{x_2}$\quad\quad For $0<\epsilon_0\ll 1$,
$P_\pm(x;z)$ is purely outgoing for any $z\in{\rm SIP}$.

\vspace{0.3cm}

Permuting $c_{11}$ and $c_{22}$ in $(C_1)'$ to get
\begin{align*}
  (C_1)\quad \gamma_2\leq 0,
\end{align*}
we obtain
\[
  {\rm (BFJ)}_{x_2}\quad\quad    (C_1)\ \&\ (C_2)\ \&\ [(C_3)_1\ |\ (C_3)_2],
\]
where $\&$ and $|$ denote the logical ``and'' and ``or'' operators,
respectively. 
\begin{remark}
  \label{rem:g2}
  Condition (BFJ)$_{x_2}$ is stronger than (BFJ)$_{x_2}^0$. See Medium (II) in
  section 3.3 with $\gamma_2>0$ but $\min\Delta(\xi^2)\geq 0$. Moreover, it is
  straightforward to deduce $\gamma_2<0$ if $\alpha_1\leq 0$.
\end{remark}
We prove that (BFJ)$_{x_2}$ is a sufficient and necessary condition of
(OC)$_{x_2}$ in the following. First, we identify the branch cuts of
$\sqrt{\mu_\pm^2}$.
\begin{lemma}
  \label{lem:bfjx2+}
For any $0<\epsilon_0\ll 1$, 
 \[
   {\rm Im}(\mu_\pm^2(z^2(\xi^2;\epsilon_0)))\geq 0,\quad\forall \xi\in\mathbb{R},
 \]
 holds if and only if condition $({\rm BFJ})_{x_2}$ holds. The equality holds
 iff $\xi=0$.
\begin{proof}
``If''. Case $\xi=0$ is trivial. Suppose $\xi>0$ in the following.
  (BFJ)$_{x_2}$ implies
  \[
    \gamma_2\leq 0 \ \&\ \alpha_2\geq 0\ \&\ \beta_1<0\ \&\ \left[\alpha_1\geq 0\ |\
      \alpha_1^2-\alpha_0\alpha_2\leq 0 \right].
  \]
  We claim that \eqref{eq:realroot} cannot hold. This is clear for $\alpha_1\geq
  0$. Suppose $\alpha_1<0$ and $\alpha_1^2\leq \alpha_0\alpha_2$ now. Thus
  $\alpha_2>0$ and
\begin{align*}
  &\alpha_1^2 - \alpha_2\left[ \frac{-\gamma_2}{c_{11}c_{33}}+\beta_1^2\epsilon^2\right]\\
  =& c_{11}^{-1}c_{33}^{-1}(c_{12}+c_{33})^2\beta_1^2[(c_{12}+c_{33})^2
  -(c_{11}-c_{33})(c_{22}-c_{33})] - \alpha_2\beta_1^2\epsilon^2<0,
\end{align*}
and again (\ref{eq:realroot}) cannot hold. Consequently, $\mu_\pm^2(z^2)$ never
touch the real axis unless $\xi= 0$. At $\xi^2=+\infty$, it is not hard to
see that
\[
  {\rm Im}(\mu_\pm^2(z^2))\eqsim \frac{-\beta_1\mp
    \sqrt{\alpha_2}}{2c_{33}c_{22}}\epsilon >0,
\]
so that the proof can be concluded by noticing that $\mu_\pm^2(z^2)$ is continuous for $z\in$SIP.

\noindent ``Only if''. We claim that $\Delta(\xi^2)$ should be nonnegative for
all $\xi\in\mathbb{R}$ first. Suppose, otherwise, for some $\xi_0\in\mathbb{R}$,
\[
\Delta(\xi_0^2)<0.
\]
Thus, one of $\pm\sqrt{\Delta(\xi_0^2)}$ must have strictly negative
imaginary part so that for $\epsilon\ll 1$,
\[
  \min({\rm Im}(\mu_\pm^2(\xi_0^2-\bi\epsilon)) <0,
\]
a contradiction. Consequently, $\alpha_2\geq 0$ so that $(C_2)$ holds, and one
of $(C_3)_1$ and $(C_3)_2$ must hold by $(C_0)$.

We only need to verify $(C_1)$. We prove it by contradiction. Suppose $(C_1)$ does not hold, i.e., $\gamma_2>0$. Then, $\beta_1^2\alpha_0 - \alpha_1^2=-4c_{22}c_{33}\gamma_2<0$,
so that 
\[
  \beta_1\sqrt{\alpha_0} +|\alpha_1|>0,
\]
as $\beta_1<0$. Thus,
\[
  \alpha_2\alpha_0 - \alpha_1^2\leq \beta_1^2\alpha_0 - \alpha_1^2<0,
\]
so that $(C_3)_1$ does not hold. Then, $(C_3)_2$ must hold so that $\alpha_1>
0$. 
If $\alpha_0>0$, then $\Delta(\xi^2)\geq \Delta(0)>0$ so that
\[
  {\rm Im}(\sqrt{\Delta(\xi^2-\bi\epsilon)}) =-\frac{\alpha_2\xi^2+\alpha_1}{\sqrt{\Delta(\xi^2)}}\epsilon + {\cal
    O}(\epsilon^2).
\]
Consider
\begin{align*}
  F(\xi^2):=&\beta_1^2\Delta(\xi^2) - (\alpha_2\xi^2+\alpha_1)^2=(\beta_1^2-\alpha_2)\alpha_2\xi^4 + 2(\beta_1^2-\alpha_2)\alpha_1\xi^2+(\beta_1^2\alpha_0-\alpha_1^2).
\end{align*}
The minimum of $F(\xi^2)$ is negative since $\beta_1^2\alpha_0<\alpha_1^2$.
There exists $\xi_*\in\mathbb{R}\backslash\{0\}$ such that
\[
  |\beta_1|<\frac{|\alpha_2\xi_*^2 + \alpha_1|}{\sqrt{\Delta(\xi_*^2)}}.
\]
It is straightforward to verify that for $\epsilon\ll1$,
\[
  {\rm Im}(\mu_{+}^2(\xi_*^2-\bi\epsilon))\leq \frac{1}{2}\left[ |\beta_1| - \frac{|\alpha_2\xi_*^2 + \alpha_1|}{\sqrt{\Delta(\xi_*^2)}}\right]\epsilon<0,
\]
a contradiction. Thus, we must have $\alpha_0=0$ so that $\xi=0$ becomes a
branch point of $\sqrt{\Delta(\xi^2-\bi\epsilon)}$. Nevertheless, for
$0<\xi^2\leq \epsilon_0$, $\epsilon = \xi^2\epsilon_0$ so that
\[
  \sqrt{\Delta(\xi^2-\bi\epsilon)} = |\xi|\sqrt{2\alpha_1(1-\epsilon_0\bi) +
    \alpha_2(1-\epsilon_0\bi)^2\xi^2} =
  |\xi|\sqrt{2\alpha_1}\sqrt{1-\epsilon_0\bi}+ {\cal O}(|\xi|^3).
\]
The proof is concluded by observing 
\[
  {\rm Im}(\mu_+^2(\xi^2-\bi\epsilon)) = -\frac{|\xi|\sqrt{2\alpha_1}\epsilon_0}{2}+ {\cal O}(\epsilon_0^{3/2})<0,
\]
\end{proof}
\end{lemma}
Lemma~\ref{lem:bfjx2+} suggests to use the positive real axis as the branch cut
of $\sqrt{\mu_\pm^2}$ so that $\mu_\pm(z^2)\in \overline{\mathbb{C}^{++}}$ for
all $z\in$SIP. Consequently, we obtain the equivalence of (BFJ)$_{x_2}$ and
(OC)$_{x_2}$.
\begin{theorem}
  \label{thm:x2}
  For $0<\epsilon_0\ll 1$, ${\rm (BFJ)}_{x_2} \Longleftrightarrow {\rm (OC)}_{x_2}$.
\end{theorem}

\subsubsection{A weak condition}
Condition (OC)$_{x_2}$ requires $\mu_\pm(z^2)\in \overline{\mathbb{C}^{++}}$
for all $z\in$SIP for $0<\epsilon_0\ll 1$, which
implies $\mu_\pm(\xi^2)\in \overline{\mathbb{C}^{++}}$ for all
$\xi\in\mathbb{R}$ by the definition (\ref{eq:mu:def}), but not vice versa.
Taking $\epsilon_0=0$ to replace SIP in (OC)$_{x_2}$ by the real axis, we obtain
a weaker condition
\begin{center}
(OC)$^0_{x_2}$:\quad\quad(OC)$_{x_2}$ with $\epsilon_0=0$.
\end{center}
To justify the equivalence of (OC)$^0_{x_2}$ and (BFJ)$_{x_2}^0$, we need a more
delicate lemma than Lemma~\ref{lem:bfjx2+}.
\begin{lemma}
  \label{lem:wkx2}
  For all $\xi\in\mathbb{R}$,
  \begin{align}
    \label{eq:weakres1}
    &\arg(\mu_-^2(z^2(\xi^2;\epsilon_0)))\in [0,\pi),\\
    \label{eq:weakres2}
    &\arg(\mu_+^2(z^2(\xi^2;\epsilon_0)))\in (-\delta_0\epsilon_0,\pi+\delta_0\epsilon_0),
 \end{align}
for $0<\epsilon_0\ll 1$, if and only if condition
 $({\rm BFJ})^0_{x_2}$ holds. Here $\delta_0>0$ is a constant independent of $\epsilon_0$.
    \begin{proof}
 ``If''. As in the proof of Lemma~\ref{lem:bfjx2+}, $(C_3)_1$ implies
      $(C_1)$ so that (\ref{eq:weakres1}) and (\ref{eq:weakres2}) follow
      immediately. We now assume $(C_3)_2$ holds but $(C_3)_1$ and $(C_1)$ not.
      Thus, $\gamma_2>0$, $\alpha_2\geq 0$, $\alpha_1\geq 0$, and $\beta_1<0$.
      By (\ref{eq:del:z2}) and by $\min_{\xi\in\mathbb{R}}\Delta(\xi^2)>0$,
      ${\rm Im}(\sqrt{\Delta(z^2)})\leq 0$ so that
      \[
        {\rm Im}(\mu_-^2(z^2))\geq \frac{-\beta_1\epsilon}{2c_{33}c_{22}}>0,
      \]
      for all $0\neq z\in$SIP. In other words, $\mu_-^2(z^2)$ never touches the
      real axis and $\arg\mu_-^2(z^2)\in (0,\pi)$. We now consider
      $\arg\mu_+^2(z^2)$. Observe that \eqref{eq:realroot} has only one positive
      root
      \[
        \xi_*^2(\epsilon) = \frac{-\alpha_1 + \sqrt{\alpha_1^2 +
            \frac{\alpha_2\gamma_2}{c_{11}c_{33}}}}{\alpha_2} + {\cal O}(\epsilon^2),
      \]
      where the prefactor in ${\cal O}$-term is independent of $\epsilon^2$ and
      the first term on the right-hand side takes a limit when $\alpha_2=0$.
      Thus, $\mu_+^2(\xi^2)$ must cross the positive real axis once only at
      $\xi_*^2$. 

      As in the proof of Lemma~\ref{lem:bfjx2+}, ${\rm Im}(\mu_+^2)>0$ at
      $\xi^2=\infty$. Thus, we can find a constant $\delta_1>0$, such that ${\rm
        Im}(\mu_+^2(z^2))>0$ for $\xi^2-\xi_*^2(0)\geq \delta_1$.
      For $\xi^2\in [0,\xi_*^2(0)+\delta_1]$, it can be easily seen
      that $|{\rm Im}(\mu_+^2(z^2))|\leq \delta_2\epsilon$ for some constant
      $\delta_2$ independent of $\epsilon_0$. If $|{\rm Re}(\mu_+^2(z^2))|$ is
      uniformly bounded below from some constant greater than $0$ w.r.t.
      $\epsilon_0$, then the proof is concluded. This indicates that we only
      need to consider the case when $\xi^2$ is near $c_{11}^{-1}$ or
      $c_{33}^{-1}$, where $\mu_+^2(\xi^2)$ is zero.

      Since $\mu_-^2(\xi^2)$ is strictly decreasing, it has exactly one root
      considering $\mu_-^2(0)\mu_-^2(+\infty)<0$ and hence $\mu_+^2(\xi^2)$ also
      has exactly one root. We thus distinguish two cases:\\
      \noindent(i) $\mu_-^2(c_{11}^{-1})=0$ such that $\mu_+^2(c_{33}^{-1})=0$.
      Thus, $\mu_+^2(c_{11}^{-1})>0$ and $\mu_-^2(c_{33}^{-1})<0$ considering
      $\Delta(c_{jj}^{-1})>0$ for $j=1,3$. But $\mu_-^2$ is decreasing so that
      $c_{11}>c_{33}$. Now suppose $\xi^2=c_{33}^{-1}+\delta$ for $|\delta|\leq
      \delta_3$, a sufficiently small constant independent of $\epsilon_0$. Then,
      \begin{align*}
        \mu_+^2(z^2) 
        =&\frac{c_{11}c_{33}\delta^2 + (c_{11}-c_{33})\delta - c_{11}c_{33}\epsilon^2+ [(c_{33}-c_{11})-2c_{11}c_{33}\delta]\epsilon\bi}{c_{33}c_{22}\mu_-^2(z^2)}.
      \end{align*}
      Here, the numerator is strictly below the real axis and
      $\arg\mu_-^2(z^2)\in(\pi - \delta_4\epsilon,\pi)$ for some constant
      $\delta_4$ independent of $\epsilon_0$.
      Thus, $\arg\mu_+^2(z^2)\in (0,\pi+\delta_4\epsilon)$.\\
\noindent(ii) $\mu_-^2(c_{33}^{-1})=0$ such that $\mu_+^2(c_{11}^{-1})=0$. One
can just permute $c_{11}$ and $c_{33}$ in the arguments above to obtain the same
conclusion. 

(\ref{eq:weakres2}) follows by combining the above results. 

``Only if''. Firstly, $\Delta(\xi^2)\geq 0$ for all $\xi\in\mathbb{R}$, since,
otherwise, one of $\mu_\pm^2(\xi_0^2)$ must lie strictly below the real axis for
some $\xi_0\in\mathbb{R}$ and (\ref{eq:weakres2}) fail at
$z^2=\xi_0^2-\bi\epsilon$. Thus, $(C_2)(i)$ and one of $(C_3)_1$ and $(C_3)_2$
must hold, implying (BFJ)$_{x_2}^0$ by Remark~\ref{rem:bfjx20}. 

     \end{proof}
\end{lemma}
The above lemma indicates that we can choose the positive real axis as the
branch cut of $\sqrt{\mu_-^2}$ and a ray $[0,\infty)e^{-\bi \theta_0}$ in
Quadrant IV as the branch cut of $\sqrt{\mu_+^2}$, where $\theta_0>0$ is a sufficiently small constant independent of $\epsilon_0$. Letting
$\epsilon_0\to 0^+$, we obtain the following theorem immediately.
\begin{theorem}
  \label{thm:x2w}
    {\rm(OC)}$^0_{x_2}\Longleftrightarrow${\rm (BFJ)}$^0_{x_2}$.
\end{theorem}

\subsection{Partially incoming Green's tensor for
  $\min_{\xi\in\mathbb{R}}\Delta(\xi^2)<0$}
By Theorem \ref{thm:x2w}, (OC)$_{x_2}^0$ does not hold only if (BFJ)$_{x_2}^0$
is violated, and hence $\min_{\xi\in\mathbb{R}}\Delta(\xi^2)<0$. We now justify
that the previously proposed two principles (P1) and (P2) are still applicable
for any general tensor $C$ satisfying $(C_0)$ but not (BFJ)$_{x_2}^0$. As in the
proof of Lemma~\ref{lem:wkx2}, $(C_3)_1$ must not hold, i.e.,
$\alpha_1^2-\alpha_0\alpha_2>0$. Then, breaking (BFJ)$_{x_2}^0$ indicates that
the following three inequalities
\[
  \beta_1<0,\quad\alpha_2\geq 0,\quad{\rm and}\quad \alpha_1\geq 0,
\]
cannot hold simultaneously. Two possible situations arise:
$\min_{\xi\in\mathbb{R}}\Delta(\xi^2)\in(-\infty,0)$ and
$\min_{\xi\in\mathbb{R}}\Delta(\xi^2)=-\infty$. The related results are
presented in the following two theorems.


  \begin{theorem}
    \label{thm:4.2.1}
    Suppose $\min_{\xi\in\mathbb{R}}\Delta(\xi^2)\in(-\infty,0)$ and
    $0<\epsilon_0\ll 1$. The positive real
    axis simultaneously serves as the branch cuts of $\sqrt{\Delta}$ and $\sqrt{\mu_\pm^2}$,
    fulfilling the two principles (P1) and (P2). The plane wave
    $P_+(x;\xi)$ is purely outgoing for all $\xi\in\mathbb{R}$, while the plane
    wave $P_-(x;\xi)$ is outgoing for all
    $\xi^2\in\mathbb{R}/(\xi_1^2(0),\xi_2^2(0))$, but incoming for $\xi^2\in
    (\xi_1^2(0),\xi_2^2(0))$, where $\xi_1^2(0)$ and $\xi_2^2(0)$ are the two positive roots of (\ref{eq:realroot}) for $\epsilon_0=0$.
    \begin{proof}
Since $\min_{\xi\in\mathbb{R}}\Delta(\xi^2)\in(-\infty,0)$, $\alpha_2>0$, $\alpha_1<0$, $\beta_1<0$ and $\gamma_2<0$ (c.f.
Remark~\ref{rem:g2}). By (\ref{eq:del:z2}),
$\Delta(z^2)$ for $z\in$SIP never crosses the positive real axis,
which shall serve as the branch cut of $\sqrt{\Delta}$. Thus, for $\xi\neq 0$,
${\rm Im}(\sqrt{\Delta(z^2)})>0$ so that ${\rm Im}(\mu_+^2(z^2))>0$ and
$\mu_+^2(z^2)$ lie completely above the real axis. We can simply choose the
positive real axis as the branch cut of $\sqrt{\mu_+^2}$.

As for $\sqrt{\mu_-^2}$, (\ref{eq:realroot}) indicates that $\mu_-^2(z^2)$ cross
the real axis twice at $\xi_j^2(\epsilon),j=1,2$ with
$\xi^2_2>-\frac{\alpha_1}{\alpha_2}>\xi_1^2$ for $0\neq z\in$SIP. Recall
$\xi_0^2$ defined in (\ref{eq:xi20}). Eq. (\ref{eq:realroot}) implies
  \begin{equation}
    \label{eq:xi0:xij}
    \Delta(\xi_0^2) - \Delta(\xi_j^2) = \frac{(c_{11}-c_{33})^2\beta_1^2}{4c_{11}^2c_{33}^2} +\beta_1^2\epsilon^2>0,\quad j=1,2,
  \end{equation}
  so that we must have $\xi_2^2>\xi_1^2>\xi_0^2$ or $\xi_0^2>\xi_2^2>\xi_1^2$.
  We claim that $\xi_0^2<\xi_1^2$. Suppose, otherwise, $\xi_0^2>\xi_2^2$. Then,
  $\xi_0^2>-\frac{\alpha_1}{\alpha_2}$ implies
  \[
    0<\alpha_2\xi_0^2+\alpha_1 = \beta_1(\beta_0+\beta_1\xi_0^2).
  \]
  Thus, $\beta_0+\beta_1\xi_0^2<0$. This together with $\alpha_1<0$ and
  $\alpha_1^2>\alpha_0\alpha_2$ implies 
  \[
    (c_{11}-c_{33})(c_{22}-c_{33})<(c_{12}+c_{33})^2<\max_{j=1,2}\left\{\frac{c_{jj}-c_{33}}{c_{jj}+c_{33}}(c_{11}c_{22}-c_{33}^2)\right\}.
  \]
  One may check that the above inequality is impossible for any $C$ satisfying $(C_0)$. 

  For $\xi_0^2< \xi_1^2$, (\ref{eq:rootformula}) implies that $\mu_-^2(z^2)$
  crosses the negative real axis twice at $\xi_j^2, j=1,2$ for $0\neq z\in$SIP.
  In other words, $\mu_-^2$ never crosses the positive real axis, which can
  serve as the branch cut of $\sqrt{\mu_-^2}$. As ${\rm Im}(\mu_-^2(z^2))\geq 0$
  for $\xi^2\ll 1$ or $\xi^2 \gg 1$, $\mu_-^2$ lies below the real axis only for
  $\xi^2\in(\xi_1^2(\epsilon),\xi_2^2(\epsilon))$. The proof is concluded by
  letting $\epsilon_0\to 0^+$.
     \end{proof}
  \end{theorem}

  \begin{theorem}
    \label{thm:4.2.2}
    Suppose $\min_{\xi\in\mathbb{R}}\Delta(\xi^2)=-\infty$ and $0<\epsilon_0\ll 1$. The negative real
    axis and the positive real axis serve as the branch cuts of
    $\sqrt{\Delta(z^2)}$ and $\sqrt{\mu_-^2(z^2)}$, respectively, for $z\in$SIP.
    The branch cut of $\sqrt{\mu_+^2(z^2)}$ should be chosen as ray $[0,+\infty)e^{-\bi
      \theta_0}$ in Quadrant IV if $\gamma_2>0$ or the positive real axis if
    $\gamma_2\leq 0$. They fulfill the two principles (P1) and (P2). The plane
    wave $P_+(x;\xi)$ is purely outgoing for all $\xi\in\mathbb{R}$, while the
    plane wave $P_-(x;\xi)$ is outgoing for $\xi^2\in[0,\xi_*^2(0)]$ but
    incoming for $\xi^2\in (\xi_*^2(0),\infty)$, where $\xi_*^2(0)$ is the
    single positive root of (\ref{eq:realroot}) for $\epsilon_0=0$.
    \begin{proof}
Since $\min_{\xi\in\mathbb{R}}\Delta(\xi^2)=-\infty$, $\alpha_2\leq 0$. We claim that $\Delta(z^2)$ never crosses the
negative real axis as $z$ travels along an SIP for $\epsilon_0\ll 1$. Otherwise,
by (\ref{eq:del:z2}), $\xi^2 = -\frac{\alpha_1}{\alpha_2}$ and
$\Delta(\xi^2)=\max_{t\in\mathbb{R}}\Delta(t^2)\leq 0$. Thus, $\alpha_0=0$
so that $\alpha_1>0$. But now $\Delta(t^2)> 0$ for $t^2\ll 1$, a
contradiction. Consequently, the negative real axis can serve as the branch cut
of $\sqrt{\Delta(z^2)}$ for $z\in$SIP.
As for the branch cuts of $\sqrt{\mu_\pm^2}$, we distinguish two cases: (i)
$\gamma_2>0$; (ii) $\gamma_2\leq 0$.

{\bf (i).} $\gamma_2>0$ so that $\alpha_1>0$, $\alpha_2<0$, and
\[
  \beta_1^2\alpha_0 - \alpha_1^2 = -4c_{22}c_{33}\gamma_2<0.
\]
Moreover, (\ref{eq:realroot}) must have two positive roots denoted by $\xi_1^2(\epsilon)$ and
$\xi_2^2(\epsilon)$. This indicates that $\mu_\pm^2(z^2)$ together cross the real axis
twice at $\xi^2=\xi_j^2, j=1,2$ with
$0<\xi_1^2<-\frac{\alpha_1}{\alpha_2}<\xi_2^2$. Thus, we must have
$\beta_1\sqrt{\alpha_0} -\alpha_1 < 0$ and $\beta_1\sqrt{\alpha_0} + \alpha_1 >
0$. For $\xi^2\ll 1$, $\epsilon = \epsilon_0|\xi|^2$. If $\alpha_0> 0$,
\begin{equation}
  \label{eq:closeto0:+}
  2c_{22}c_{33}\mu_\pm(z^2) \eqsim \beta_0 -\beta_1\bi\epsilon \pm
  \sqrt{\alpha_0 -2\alpha_1\epsilon\bi} \eqsim \beta_0\pm
  \sqrt{\alpha_0} -(\beta_1\pm \frac{\alpha_1}{\sqrt{\alpha_0}})\bi\epsilon.
\end{equation}
If $\alpha_0=0$,
\begin{equation}
  \label{eq:closeto0:0}
2c_{22}c_{33}\mu_\pm(z^2) \eqsim \beta_0 -\beta_1\bi\epsilon \pm
|\xi|\sqrt{2\alpha_1}\sqrt{1-\epsilon_0\bi}\eqsim \beta_0 \mp \frac{|\xi|\sqrt{2\alpha_1}}{2}\epsilon_0\bi.
\end{equation}
At $\xi^2=+\infty$,
\begin{equation}
  \label{eq:infty}
  2c_{22}c_{33}\mu_\pm(z^2) \eqsim \beta_1\xi^2\pm \sqrt{|\alpha_2|}\xi^2\bi.
\end{equation}
Therefore, either of $\mu_+^2$ and $\mu_-^2$ crosses the real axis exactly
once. Suppose $\beta_1\neq 0$ for the moment. By
\[
  2c_{33}c_{22}{\rm Im}(\mu_\pm^2(-\alpha_1/\alpha_2 -\bi\epsilon)) =
  -\beta_1\epsilon,
\]
it must hold that $\mu_+^2(\xi_2^2-\bi\epsilon)$ and
$\mu_-^2(\xi_1^2-\bi\epsilon)\in\mathbb{R}$ if $\beta_1>0$, or that
$\mu_+^2(\xi_1^2-\bi\epsilon)$ and $\mu_-^2(\xi_2^2-\bi\epsilon)\in\mathbb{R}$
if $\beta_1<0$. On the other hand, (\ref{eq:xi0:xij}) and 
\[
  \Delta(\xi_1^2) = \Delta(\xi_2^2)=\frac{\gamma_2}{c_{11}c_{33}} - \beta_1^2\epsilon^2 +
  \alpha_0 >0, 
\]
imply that $\xi_1^2<\xi_0^2<\xi_2^2$. Consequently, (\ref{eq:rootformula})
implies $\mu_-^2(z^2)$ crosses the negative real axis once while $\mu_+^2(z^2)$
crosses the positive real axis once. We can use the positive real axis as the
branch cut of $\sqrt{\mu_-^2}$. As for $\sqrt{\mu_+^2}$, it holds that
$2c_{33}c_{22}{\rm Re}(\mu_+^2)\geq \beta_0$ and $|{\rm Im}(\mu_+^2)|\leq
\delta_2\epsilon$ for $\xi^2\in[0,\xi_2^2+\delta_1]$, where constants $\delta_1$
and $\delta_2$ are independent of $\epsilon_0$. Thus, we can use ray
$[0,\infty)e^{-\bi\theta_0}$ in Quadrant IV as the branch cut of
$\sqrt{\mu_+^2}$. If $\beta_1=0$, since ${\rm Re}(\sqrt{\Delta})\geq 0$ for all
$z\in$SIP, $\mu_+^2(\xi_0^2-\bi\epsilon)>0$ so that
$\mu_-^2(\xi_0^2-\bi\epsilon)<0$ by Lemma~\ref{lem:realcros}. The rest is
similar as before.

{\bf (ii).} $\gamma_2\leq 0$ so that $\beta_1<0$. Suppose first
$\alpha_2<0$. Then, $\mu_\pm^2(z^2)$ together cross the real axis once at
$\xi_2^2$ with $\xi^2_2>\max(0, -\frac{\alpha_1}{\alpha_2})$ for $0\neq z\in$SIP
by (\ref{eq:realroot}). But, (\ref{eq:xi0:xij}) implies $\xi_0^2<\xi_2^2$ so
that the crossing point must occur on the negative real axis so that none of
$\mu_\pm^2(z^2)$ cross the positive real axis for $0\neq z\in$SIP. Thus, the
positive real axis can serve as the branch cuts of $\sqrt{\mu_\pm^2}$. If $\alpha_2=0$ so that $\alpha_1<0$ and $\alpha_0>0$. Again,
(\ref{eq:realroot}) has only one positive root $\xi_2^2$ with $\xi_2^2>\xi_0^2$. The rest is similar as before.

To sum up, the above theory indicates that $\mu_-^2(z^2)$ lies below the real
axis only when $\xi^2>\xi_*^2(\epsilon)$, where
\begin{equation}
  \xi_*^2(\epsilon) = \begin{cases}
   \xi_1^2(\epsilon) & \beta_1>0, \gamma_2>0;\\  
   \xi_2^2(\epsilon) & \beta_1<0;\\ 
   \xi_0^2(\epsilon) & \beta_1=0, \gamma_2>0.
    \end{cases}
\end{equation}
The proof is concluded by letting $\epsilon_0\to 0^+$.
     \end{proof}
  \end{theorem}
All previous results are summarized in Table~\ref{tab:bc}. 
\begin{table}
\centering
\begin{tabular}{|l|c|c|c|c|} 
\hline
  $\min_{\xi\in\mathbb{R}}\Delta(\xi^2)$                  & \multicolumn{2}{c|}{$(-\infty,0)$}                  & $[0,+\infty) $(Outgoing) & $\{-\infty\}$                                   \\ \hline
B.C. $\sqrt{\Delta(z^2)}$                  & \multicolumn{2}{c|}{$[0,+\infty)$}                  & \multicolumn{2}{c|}{$(-\infty,0)$/entire}                   \\ 
\hline
B.C. $\sqrt{\mu_-^2(z^2)}$ & \multicolumn{4}{c|}{$[0,+\infty)$}                                    \\ 
\hline
\multirow{2}{*}{B.C. $\sqrt{\mu_+^2(z^2)}$} & \multirow{2}{*}{$[0,+\infty)$} & \multicolumn{3}{l|}{$\gamma_2\leq 0$: $[0,+\infty)$} \\  
  \cline{4-5} &                     & \multicolumn{3}{l|}{$\gamma_2>0$: $[0,+\infty)e^{-\bi\theta_0}$} \\
\hline
\end{tabular}
\caption{Branch cuts of the three square-root functions $\sqrt{\Delta(z^2)}$,
  $\sqrt{\mu_-^2(z^2)}$ and $\sqrt{\mu_+^2(z^2)}$ for $z\in$SIP. Here, B.C.
  stands for branch cut and
  $\theta_0>0$ is a sufficiently small constant independent of $\epsilon_0$; $P_\pm(x;\xi)$ is outgoing along $x_2$-direction only when
  $\min_{\xi\in\mathbb{R}}\Delta(\xi^2)\in[0,+\infty)$; $\sqrt{\Delta(z^2)}$ can
  be entire if $\min_{\xi\in\mathbb{R}}\Delta(\xi^2)=0$.}
\label{tab:bc}
\end{table}
To conclude this subsection, we make the following useful remark. 
\begin{remark}
  \label{rem:gensip}
  Let $a(\xi)$ be a continuous function for $\xi\in\mathbb{R}$ such that ${\rm
    sgn}(a(\xi)) = -{\rm sgn}(\xi)$ (hence $a(0)=0$), and
  $0<\sup_{\xi\in\mathbb{R}} |a(\xi)\xi^{-1}|\ll 1$.
  It is not hard to derive that the previously used {\rm SIP} can be replaced by
  the more general {\rm SIP}
    $\{z = \xi + a(\xi)\bi: \xi\in\mathbb{R}\}$,  as it ensures ${\rm Im}(z^2)= 2\xi a(\xi) \ll \xi^2 - a^2(\xi)={\rm Re}(z^2)$,
  which all previous arguments only require. 
\end{remark}

\subsection{Green's tensor in dissipative medium and limiting absorption principle}
We are ready to justify the physical correctness of the prescribed Green tensor
$G$ by the limiting absorption principle. We shall overload $\mu_\pm^2(z^2)$,
$\Delta(z^2)$, and $G(x-y;\rho\omega^2)$ with $\mu_\pm^2(z^2;\rho\omega^2)$,
$\Delta(z^2;\rho\omega^2)$, and $G(x-y)$ respectively to emphasize their
dependence on $\rho\omega^2$ in this subsection. Clearly,
\begin{equation}
  \label{eq:rel}
  \mu_\pm^2(z^2;\rho\omega^2) = \rho\omega^2\mu_\pm^2(z^2/(\rho\omega^2);1),\quad
  \Delta(z^2;\rho\omega^2) = \rho^2\omega^4\Delta(z^2/(\rho\omega^2);1).
\end{equation}
For simplicity, we suppress the argument $\rho\omega^2$ if it equals $1$ in the following. Now, let
$\rho\omega^2 = 1 + \bi \delta$ for $0< \delta\ll 1$. Then, by (\ref{eq:rel})
and by Remark~\ref{rem:gensip},
we may compute the related Green tensor $G$ based on the same idea in section 2
but directly use the real path, considering that $\{\xi(1+\bi\delta)^{-1/2}: \xi\in\mathbb{R}\}$ is
a suitable SIP to well define $\sqrt{\mu^2_\pm(\xi^2/(1+\bi\delta))}$ and
$\sqrt{\Delta(\xi^2/(1+\bi\delta))}$ according to Theorems~\ref{thm:x2w}, \ref{thm:4.2.1}, and
\ref{thm:4.2.2}. Thus, we can define for $\xi\in\mathbb{R}$ that
\begin{equation}
  \label{eq:def:1}
\mu_\pm(\xi^2;1+\bi\delta) := (1+\bi\delta)^{1/2}\sqrt{\mu^2_\pm(\xi^2/(1+\bi\delta))},
\end{equation}
based on the prescribed branch cuts in Table~\ref{tab:bc}. Nevertheless, it
turns out that the branch cuts of $\sqrt{\mu^2_\pm(\xi^2;1+\bi\delta)}$ are much
simpler than $\sqrt{\mu^2_\pm(\xi^2/(1+\bi\delta))}$ due to the following lemma.
\begin{lemma}
  \label{lem:dist}
For any $\delta>0$, the pathes $P_\pm^\delta:=\{\mu_\pm^2(\xi^2;1+\bi\delta):\xi\in\mathbb{R}\}$ never touch the positive real
axis. Thus, the following distance function
\begin{equation}
  \label{eq:dist}
  {\rm dist}(P_\pm^\delta,\mathbb{R}^+)\geq c(\delta),\quad \mathbb{R}^+=[0,\infty),
\end{equation}
for some strictly positive function $c(\delta)$ that approaches
$0$ as $\delta\to 0^+$.
\begin{proof}
  Observe (\ref{eq:mupm}). At most one of $\mu^2_\pm$ touches the positive real axis at $\xi^2=\xi_0^2\geq
  0$ and it must hold that
  \begin{equation}
    \label{eq:tmp3}
    \Delta(\xi_0^2;1+\bi\delta) = (a - \bi\delta \beta_0)^2,\quad
    \beta_0+\xi_0^2\beta_1 + a > 0.
  \end{equation}
  for some costant $a$ depending on $\xi_0^2$. Compare the imaginary part to
  obtain $a = -(\alpha_0 + \alpha_1\xi_0^2)\beta^{-1}_0$, and then the real part
  to obtain that $\xi_0^2$ must be a root of
  \[
    F(t):=(\alpha_2 - \alpha_1^2\beta_0^{-2})t^2 +
    2\alpha_1\beta_0^{-2}(\beta_0^2-\alpha_0)t +
    (\alpha_0\beta_0^{-2}+\delta^2)(\beta_0^2-\alpha_0) = 0.
  \]
  The inequality in (\ref{eq:tmp3}) implies $0<\xi_0^2<2(c_{11}+c_{33})^{-1}$. It
  is straightforward to verify from $(C_0)$ that
    \[
\alpha_2 - \alpha_1^2\beta_0^{-2} < 0, \beta_0^2 - \alpha_0> 0,
    \]
    Besides, it can be seen that $F(0)>0$ and
    \[
 F(2(c_{11}+c_{33})^{-1})=\frac{4c_{22}c_{33}(c_{11}-c_{33})^2}{(c_{11}+c_{33})^2}
 + 4c_{22}c_{33}\delta^2>0,
    \]
    implying $F(\xi_0^2)\geq 4c_{22}c_{33}\delta^2$ since $F$ opens
    downward. The proof is finished by noticing that
    \[
      |{\rm Im}(\mu_\pm^2(\xi^2;1+\bi\delta))|\geq \frac{(\beta_0-\sqrt{\alpha_0})\delta}{4c_{33}c_{22}}, 
    \]
    for $\xi^2>C$, a sufficiently large constant independent of $\delta$.
\end{proof}
\end{lemma}
Lemma~\ref{lem:dist} indicates that we can define $\mu_\pm(\xi^2;1+\bi\delta)$
by directly choosing the positive real axis as the branch cut of
$\sqrt{\mu_\pm^2(\xi^2;1+\bi\delta)}$. This new definition coincides with
(\ref{eq:def:1}) since both produce the same results at $\xi=\pm\infty$ and both
are analytic for $\xi\in\mathbb{R}$. Due to the strictly positive distance of
$P_\pm^\delta$ and $\mathbb{R}^+$ by (\ref{eq:dist}), the real path can
slightly move to $R_{\delta}:=\{\xi+p(\delta)\bi:\xi\in\mathbb{R}\}$ for some
positive constant $p(\delta)$ that approaches $0$ as $\delta\to 0^+$ so that
$\mu^2_\pm((\xi+p(\delta)\bi)^2;1+\bi\delta))$ are still strictly away from
$\mathbb{R}^+$ and that $\min_{\xi\in\mathbb{R}}{\rm
  Im}(\mu_\pm(\xi+p(\delta)\bi;1+\bi\delta))\geq q(\delta)$ for some positive
constant $q(\delta)$ that approaches $0$ as $\delta\to 0^+$. Consequently, by
Cauchy's theorem, the real path used in (\ref{eq:green11}-\ref{eq:green22}) can
be deformed to the new path $R_{\delta}$. Since
\[
  \max_{z\in R_{\delta}}{\rm Im}({\rm \mu}_\pm(z;1+\bi\delta) |x_2| + z
  |x_1|)\geq p(\delta) |x_1| + q(\delta)|x_2|,
\]
it is not hard to derive that for $|x|\gg 1$,
\[
  |G_{ij}(x;1+\bi\delta)||\leq C(\delta) e^{-p(\delta) |x_1|/2  - q(\delta)|x_2|/2},\quad i,j=1,2,
\]
for some positive constant $C(\delta)$ independent of $x$. Therefore, $G(x;1+\bi\delta)$ is dissipative
in the sense that it decays exponentially at $|x|=\infty$. Moreover, by (\ref{eq:rel}), (\ref{eq:def:1})
and by change of variables, it is easy to deduce that  
\[
  G(x;1) = \lim_{\delta\to 0^+} G(x;1+\bi\delta),\quad |x|\neq 0.
\]
Consequently, our proposed Green tensor for $\rho\omega^2>0$  singled out by the two fundamental
principles (P1) and (P2) is the limit of dissipative Green's tensor, thus
justifying its physical correctness. 

\subsection{Representative examples}

\begin{figure}[!ht]
  \centering
(a)\includegraphics[width=0.3\textwidth]{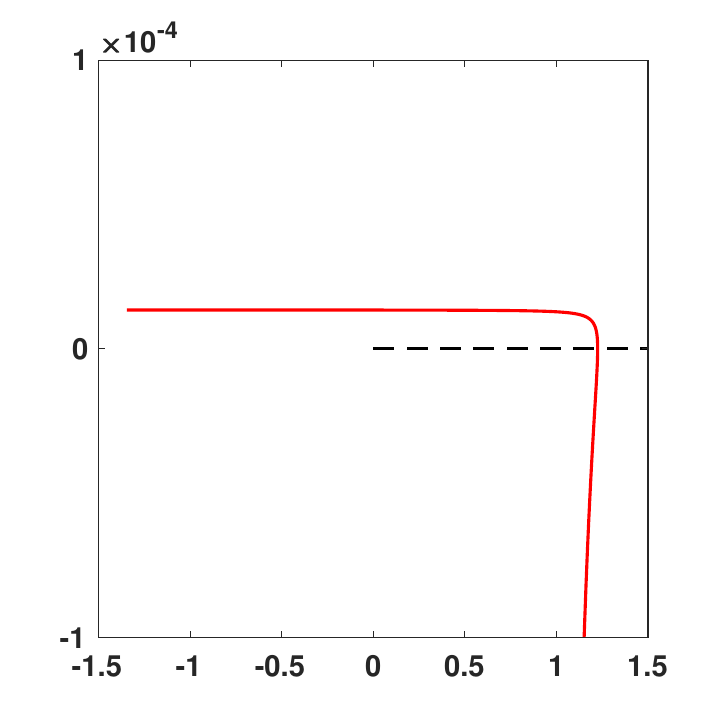}
(b)\includegraphics[width=0.3\textwidth]{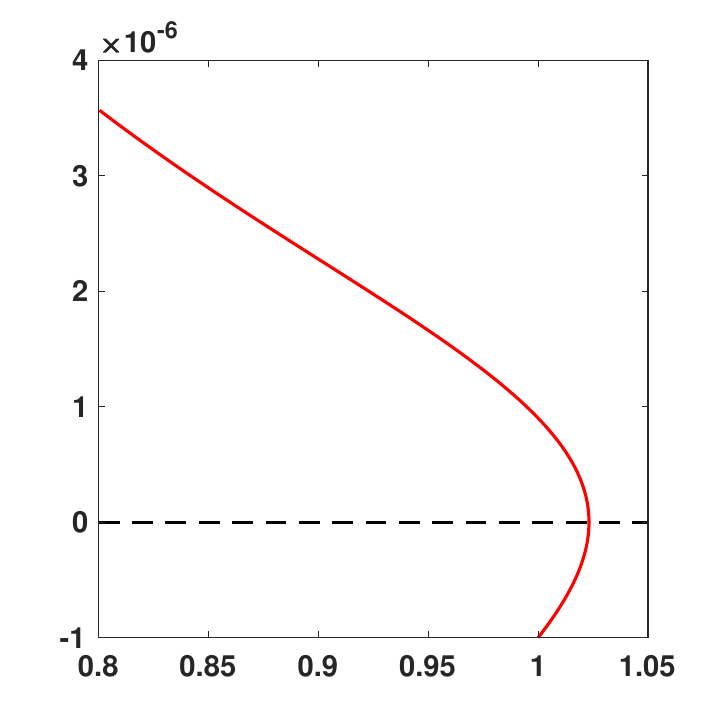}
\caption{Contours of $\mu_+^2(z^2)$ for $\epsilon_0 = 10^{-5}$ around the
  positive real axis: (a) medium (II); (b) medium (VI). Dashed lines indicate
  the positive real axis.}
  \label{fig:bc}
\end{figure}

Here, we list several representative examples below and
shall study wave scattering problems in such media in sections 5 and 6. Note
that we have excluded the trivial isotropic case.
\begin{itemize}
\item[(i).] Media with $\min\Delta(\xi^2)\geq 0$:
  \[
   {\rm (I)}\quad\quad\quad\quad c_{11} =  4, c_{22} = 20, c_{12} = 3.8, c_{33}
   = 2,
  \]
  with $\gamma_2\leq 0$;
  \[
    {\rm (II)}\quad\quad\quad\quad c_{11} = 36, c_{22} = 1, c_{12} = 2, c_{33} =
    1,
  \]
  with $\gamma_2>0$.
\item[(ii).] Medium with $\min_{\xi\in\mathbb{R}}\Delta(\xi^2)\in(-\infty,0)$: 
    \[
     {\rm (III)}\quad\quad\quad\quad  c_{11} = 0.9, c_{22} = 10, c_{12} = 0.1,
     c_{33} = 1.
    \]
\item[(iii).] Media with $\min_{\xi\in\mathbb{R}}\Delta(\xi^2)=-\infty$: 
  \begin{align*}
     &{\rm (IV)}\quad\quad\quad\quad  c_{11} = 0.9, c_{22} = 10, c_{12} = 2.5, c_{33} = 1,\quad (\beta_1>0)\\
     &{\rm (V)}\quad\quad\quad\quad c_{11} = 0.9, c_{22} = 10, c_{12} = \sqrt{10}-1, c_{33} = 1,\quad (\beta_1=0)\\
     &{\rm (VI)}\quad\quad\quad c_{11} = 0.9, c_{22} = 10, c_{12} = 2, c_{33} = 1,\quad (\beta_1<0)
  \end{align*}
  with $\gamma_2>0$;
  \begin{align*}
     &{\rm (VII)}\quad\quad\quad c_{11} = 0.9, c_{22} = 10, c_{12} = 1.835, c_{33} = 1,\quad (\alpha_1>0,\alpha_2<0)\\
     &{\rm (VIII)}\quad\quad\quad\quad c_{11} = 0.9, c_{22} = 10, c_{12} = 1.5, c_{33} = 1,\quad (\alpha_1<0,\alpha_2<0)\\
     &{\rm (IX)}\quad\quad\quad\quad c_{11} = 0.9, c_{22} = 10, c_{12} = 1, c_{33} = 1,\quad (\alpha_2=0)
  \end{align*}
  with $\gamma_2\leq 0$.
\end{itemize}
Among them, we show the contours of $\mu_+^2(z^2)$ for media (II) and (VI) with  $\gamma_2>0$, in Figure~\ref{fig:bc}. Clearly, they cross the
positive real axis, which cannot serve as their branch cuts.

\subsection{Proof of the conjecture in \cite{becfaujol03}}
In \cite{becfaujol03}, the authors proved that (BFJ)$_{x_2}^0$ is sufficient
to stabilize a PML along the $x_2$-direction for the associated time-domain
elastic wave equation. They further conjectured that
(BFJ)$_{x_2}^0$ should be necessary, the proof of which is yet incomplete. Here,
we give a rigorous proof for this conjecture, thanks to the complete 
classification of the plane waves $P_\pm(x;\xi)$ constituting Green's tensor
$G$.

In doing so, we borrow some
important notation from \cite{becfaujol03} as listed below. As suggested
therein, the stability of the time-domain PML is closely related to the roots of the following function
\begin{align}
  \tilde{F}_{\rm PML}(\omega; k, \zeta) := {\rm det}(\Gamma(\omega (\omega+\bi \zeta) k_1, \omega k_2) - \rho\omega^2(\omega+\bi\zeta)^2 I_2),
\end{align}
where $k=(k_1,k_2)\in\mathbb{R}^2$, $\zeta>0$ is a constant damping coefficient
in the time-domain PML formulation (27) of \cite{becfaujol03}, and
\begin{align}
  \Gamma(k_1,k_2) = \left( \begin{matrix}
    c_{11}k_1^2 + c_{33}k_2^2 & (c_{12}+c_{33})k_1k_2 \\
    (c_{12}+c_{33})k_1k_2 & c_{33}k_1^2 + c_{22}k_2^2
    \end{matrix} \right).
\end{align}
Note that we have revised the original term $\omega-\bi\zeta$ to
$\omega+\bi\zeta$ since we preassumed the time $e^{-\bi\omega t}$ here.
    For simplicity, we assume $\rho=1$ in the following. It is
    straightforward to derive from $\tilde{F}_{\rm PML}=0$ that $\omega$ solves
    \begin{equation}
      \label{eq:timecha}
      c_{33}c_{22}\left[ \frac{k_2^2}{(\omega+\bi\zeta)^2} \right]^2 -\left[ \beta_0 + \frac{k_1^2}{\omega^2}\beta_1 \right]\frac{k_2^2}{(\omega+\bi\zeta)^2}  + \left[ 1 - c_{11}\frac{k_1^2}{\omega^2} \right]\left[ 1 - c_{33}\frac{k_1^2}{\omega^2} \right] = 0,
    \end{equation}
    where we recall $\beta_0$ and $\beta_1$ are defined in (\ref{eq:b0}) and
    (\ref{eq:b1}) and we have presumed $\omega\notin \{0,-\bi\zeta\}$. Comparing it with (\ref{eq:cha:mu}), we see that
    $\frac{k_1^2}{\omega^2}$ and $\frac{k_2^2}{(\omega+\bi\zeta)^2}$ play the
    roles of $\frac{\xi^2}{\omega^2}$ and $\frac{\eta}{\omega^2}$,
    respectively. We have the following theorem.
\begin{theorem}
  For any fixed $\zeta>0$,  the
  roots of $\tilde{F}_{\rm PML}(\omega;k,\xi)$ satisfy ${\rm Im}(\omega)\leq
  0$ for all $k\in\mathbb{R}^2$ if and only if (BFJ)$^0_{x_2}$ holds.  
  \begin{proof}
    \noindent``Sufficient'': Though the authors in \cite{becfaujol03} have justfied the
    sufficiency of (BFJ)$^0_{x_2}$, we here give an alternative proof. Suppose
    otherwise there exists $k\in\mathbb{R}^2$ such that a root of $F_{\rm
      PML}$ satisfies ${\rm Im}(\omega)>0$. Then, there must exist
    $\tk=(\tk_1,\tk_2)\in\mathbb{R}^2$ such that $F_{\rm PML}$ has a real root
    $\tomg$. Then, (\ref{eq:timecha}) suggests that
    $\frac{\eta}{\omega^2}=\frac{\tk_2^2}{(\tomg+\bi\zeta)^2}$ satisfies 
    (\ref{eq:cha:mu}) for $\frac{\xi^2}{\omega^2} =
    \frac{\tk_1^2}{\tomg^2}$. This indicates one of $\mu_\pm^2(\xi^2)$ lies
    below the real axis for $\xi^2 = \frac{\tk_1^2\omega^2}{\tomg^2}$, a
    contradiction with Thereom~\ref{thm:x2}.
    
    ``Necessary'': Suppose now (BFJ)$_{x_2}^0$ does not hold. Then, the results
    in section 3.2 imply that for some $\xi^2 = a - \bi\delta\in$SIP with
    $a>0$ and $0<\delta\ll 1$, there exists $\zeta = b - \bi c$ with
    $b\in\mathbb{R}$ and $c\geq c_0>0$ such that
    \[
      c_{33}c_{22}\left[ b-\bi c\right]^2 -\left[ \beta_0 + (a
        -\bi\delta)\beta_1 \right](b-\bi c) + \left[ 1 - c_{11}(a-\bi\delta) \right]\left[ 1 - c_{33}(a-\bi\delta) \right] = 0.
\]
In the above, we remark that $c_0$ is independent of $\delta$ since $\zeta=b-\bi c$
must correspond to a real backward wave $P_m(x;\xi)$ as $\delta\to 0^+$. This
indicates that 
\[
  \frac{k_1^2}{\omega^2} = a - \bi \delta\quad {\rm and}\quad
  \frac{k_2^2}{(\omega+\bi\zeta)^2} = b - \bi c
\]
satisfy (\ref{eq:timecha}). We now claim that the above two equations with
three unknowns $k_1\in\mathbb{R}, k_2\in\mathbb{R}, \omega\in\mathbb{C}$ have a
solution with ${\rm Im}(\omega)>0$. Clearly, 
\begin{align*}
  \frac{k_1}{\omega} =&\sqrt{a-\bi\delta}:= a_1 - \bi \delta_1,\\
  \frac{k_2}{\omega+\bi\zeta} =& \sqrt{b-\bi c}:=b_1 - \bi c_1.
\end{align*}
for some positive constants $a_1,b_1,c_1$ and $0<\delta_1\ll 1$. Then, one finds
\[
  \frac{k_1}{a_1-\bi\delta_1} + \bi \zeta = \frac{k_2}{b_1-\bi c_1},
\]
such that
\begin{align*}
  k_1 = \frac{b_1 a_1^2 + b_1 \delta_1^2}{a_1c_1 - b_1\delta_1}\zeta>0,\quad
  k_2 = \frac{a_1(b_1^2 + c_1^2)}{a_1c_1 - b_1\delta_1}\zeta>0.
\end{align*}
Consequently,
\[
{\rm Im}(\omega) = {\rm Im}\left[  \frac{k_1}{a_1-\bi\delta_1}\right] =
\frac{k_1\delta_1}{a_1^2+\delta_1^2} > 0, 
\]
a contradiction.
\end{proof}
\end{theorem}

According to \cite[Eq. (49)]{becfaujol03}, this affirmatively resolves the
conjecture. Alternatively, from the physical perspective, Theorem~\ref{thm:x2w}
shows that (BFJ)$_{x_2}^0$ is a sufficient and necessary condition for
$P_\pm(x;\xi)$ to be purely outgoing. It is clear that such a property holds for
any $\rho>0$ and $\omega\in\mathbb{R}$. By the inverse Fourier transform w.r.t.
$\omega$, the associated time-domain Green tensor consists of purely outgoing
waves as well. According to section 1.2, this is an equivalent condition for a
PML to stably absorbing such time-domain Green's tensor and any hence the
solution of any associated Cauchy problem.

\section{Fast evaluation of Green's tensor}
In this section, we propose a fast algorithm for evaluating Green's tensor and
its gradient. Without loss of generality, we assume $x_1,x_2>0$ and shall study
$G_{11}(x)$ in detail only. Observe the coefficients $p_1^{x_1}$ and $q_1^{x_1}$
in (\ref{eq:p1}) and (\ref{eq:q1}). It can be seen that as
$\xi^2\to+\infty$, 
\begin{align*}
  p_1^{x_1}, q_1^{x_1}=\left\{
  \begin{array}{ll}
    {\cal O}(|\xi|^{-1}), & {\rm if}\quad\alpha_2\neq 0;\\
    {\cal O}(1), & {\rm if}\quad\alpha_1\neq 0, \alpha_2=0;\\
    {\cal O}(|\xi|), & {\rm if}\quad\alpha_1= 0, \alpha_2=0.
  \end{array}
  \right.
\end{align*}
This suggests to study the three cases separately. 

\noindent{\bf Case (i).} Suppose $\alpha_2=\alpha_1=0$. Then, it can be derived
based on $(C_0)$ that
\[
  0<c_{12}+2c_{33}=c_{11}=c_{22}>c_{33},
\]
Let $c_{33}=\mu>0$ and $c_{12}=\lambda$ so that $c_{11}=c_{22}=\lambda+2\mu$
with $\lambda+\mu=c_{12}+c_{33} = c_{22}-c_{33}>0$. Then, $\mu_+^2 = k_s^2 -
\xi^2$ and $\mu_-^2= k_p^2 - \xi^2$, where $k_s = \sqrt{\frac{1}{\mu}}$ and
$k_p= \frac{1}{\sqrt{\lambda+2\mu}}$. Thus,
\[
  p_1^{x_1}=\frac{\bi\mu_+}{2},\quad q_1^{x_1} = -\frac{\bi\mu_-}{2} +
  \frac{k_p^2\bi}{2}\mu_-^{-1},\quad p_2^{x_2} = -\frac{\bi \xi}{2}.
\]
By (\ref{eq:phik}), we obtain
\begin{align}
  G_{11}(x) 
  =& k_s^2\Phi_{k_s}(x) + \partial^2_{x_1}(\Phi_{k_s}(x) - \Phi_{k_p}(x)).
\end{align}
One similarly obtains
\begin{align}
  G_{21}(x)=&G_{12}(x) = \partial_{x_1x_2}^2(\Phi_{k_s}(x) - \Phi_{k_p}(x)),\\
  G_{22}(x) =& k_s^2\Phi_{k_s}(x) + \partial^2_{x_2}(\Phi_{k_s}(x) - \Phi_{k_p}(x)).
\end{align}
We have reproduced Green's tensor for an isotropic medium with Lam\'{e}
constants $\lambda$ and $\mu$. As $|x|\to 0$, $\Phi_k(x) =
-\frac{1}{2\pi}\log|x| + {\cal O}(|x|^2\log|x|)$, it is easy to conclude that
$G(x)$ exhibits logarithmic singularity as $|x|\to 0$. The evaluation of $G$ and
its gradient is straightforward.

\noindent{\bf Case (ii).} Suppose now $\alpha_1\neq 0$ and $\alpha_2=0$ so that
$\beta_1<0$. To simplify the presentation, we consider $\alpha_1>0$ only. 
We find Taylor's expansions
of $p^{x_1}_1$ and $q_1^{x_1}$ at $|\xi|=+\infty$, 
\begin{align}
  \label{eq:exppx1}
p^{x_1}_{1}(\xi^2)e^{\bi\mu_+(\xi^2) x_2}
  &=e^{\bi|\xi| p_0^+ x_2}\sum_{n=0}^{\infty}\frac{c_n^{+}(x_2) }{|\xi|^n},\\
  \label{eq:expqx1}
q^{x_1}_{1}(\xi^2)e^{\bi\mu_-(\xi^2) x_2}
  &=e^{\bi|\xi| p_0^- x_2}\sum_{n=0}^{\infty}\frac{c_n^{-}(x_2) }{|\xi|^n},
\end{align}
where $p_0^\pm = p_0 = \lim_{\xi\to\infty}\mu_\pm(\xi^2)/|\xi|=\bi \sqrt[4]{c_{11}/c_{22}}$, 
\begin{align*}
  c_0^+(x_2) =& a_0 e^{\bi b_0 x_2},\quad c_0^-(x_2) = -a_0 e^{-\bi b_0 x_2},\\
  a_0 =& \lim_{\xi\to\infty} p_1^{x_1}(\xi^2),\quad b_0 = \lim_{\xi\to\infty} (\mu_+(\xi^2) - p_0|\xi|), 
\end{align*}
and the other coefficients $c_n^\pm(x_2)$ are smooth at $x_2=0$. Then, it can be seen that
\begin{align}
  \label{eq:G11:case2}
  G_{11}(x) 
  =& \frac{1}{2\pi}\int_{-M}^{M}\left[ p^{x_1}_{1}(\xi^2)e^{\bi\mu_+(\xi^2) x_2} + q^{x_1}_{1}(\xi^2)e^{\bi\mu_-(\xi^2) x_2} \right]e^{\bi \xi x_1} d\xi\nonumber\\
  &+ \frac{1}{2\pi}\sum_{n=0}^{\infty}\left[c_n^{+}(x_2)I_M^n(x; p_0^+) + c_n^-(x_2)I_M^n(x; p_0^-)  \right],
  \end{align}
  where $M>0$ is a sufficiently large constant, and we have defined
\begin{equation}
  I_M^n(x;\alpha):=\int_{M}^{+\infty} + \int^{-M}_{-\infty}\frac{e^{\bi |\xi| \alpha x_2 + \bi \xi x_1}}{|\xi|^n} d\xi, \quad {\rm Im}(\alpha)> 0.
\end{equation}
It can be seen that, for all $n\in\mathbb{N}$,
\begin{align*}
  I_{M}^n(x;\alpha) 
  =& M^{1-n}[E_{n}(-M\bi(\alpha x_2+x_1))+ E_{n}(-M\bi(\alpha x_2-x_1))], 
\end{align*}
where $E_n(z)=z^{n-1}\int_{z}^{+\infty}t^{-n} e^{-t}dt$ is the generalized
exponential integral (c.f. (8.19.2) in \cite{nist10}). By (8.19.8) in
\cite{nist10},
\[
E_0(z) = z^{-1} + {\rm ANA},\quad E_n(z)=-\frac{(-z)^{n-1}}{(n-1)!}\log z +
{\rm ANA},\quad n\geq 1, 
\]
where ANA refers to some analytic function and the branch cut of $\log$ is
chosen as the negative real axis. Thus,
\begin{align}
  \label{eq:IM0}
  I^0_{M}(x;\alpha)  =& \frac{2\bi\alpha x_2}{\alpha^2x_2^2 - x_1^2} + {\rm ANA},\\
  \label{eq:IMn}
  I^n_{M}(x;\alpha) =& \frac{-\bi^{n-1}}{(n-1)!}\Big[(\alpha x_2 + x_1)^{n-1}\log(\alpha x_2 + x_1) 
  + (\alpha x_2-x_1)^{n-1}\log(\alpha x_2 - x_1) \Big] + {\rm ANA},
\end{align}
so that, as $|x|\to 0^+$,
\begin{equation}
  G_{11}(x) = \frac{2\bi p_0x_2a_0\sin(b_0x_2)}{\pi(p_0^2x_2^2-x_1^2)}
   - \frac{c_1^+(0)+c_1^-(0)}{2\pi}\log(-p_0^2x_2^2+x_1^2) + o(1). 
\end{equation}
Considering that $p_0=\bi|p_0|$,
\[
\left|\frac{2\bi p_0x_2a_0\sin(b_0x_2)}{\pi(p_0^2x_2^2-x_1^2)}\right|\leq \frac{2|a_0b_0|}{\pi |p_0|},
\]
so that $G_{11}$ exhibits logarithmic singularity at $z=0$. In fact, it can be
seen from (\ref{eq:G11:case2}), (\ref{eq:IM0}) and (\ref{eq:IMn}) that
\[
  G_{11}(x) = \frac{2\bi p_0x_2a_0\sin(b_0x_2)}{\pi(p_0^2x_2^2-x_1^2)} +
  {\rm ANA}\times \log(p_0 x_2 + x_1) + {\rm ANA}\times \log(p_0 x_2 -
  x_1) + {\rm ANA}.
\]
One similarly finds logarithmic singularities of $G_{12}(x)$ and $G_{22}(x)$ at $x=0$.

To numerically evaluate $G_{11}$, we directly truncate the series in
(\ref{eq:G11:case2}) at $n=N_0$. Empirically, $N_0$ and $M$ are chosen such that
$M^{-(N_0+1)}\leq \epsilon_{T}$, for some small error threshold $\epsilon_T>0$.
The functions $c_n^\pm(x_2), n=0,\cdots, N_0$ can be precomputed via symbolic
evaluation only once before the computation. The integral in
(\ref{eq:G11:case2}) can be computed via the adaptive $16$-point Gauss-Legendre
quadrature rule. If we require $G$ at a large amount of grid points of $x$, it
becomes extremely costly to directly evaluate such an integral the same amount
of times. To tackle this
difficulty, we evaluate the integral at Chebyshev points of a sufficiently large
rectangular domain enclosing the whole grid points, and then retrieve the
targeted values by the simple Lagrange interpolation, considering that
(\ref{eq:G11:case2}) is an analytic function of $x$.

To approximate the gradient of $G_{11}$, one directly differentiates
(\ref{eq:green11}) and then uses the same strategy as the above. The evaluation
of $G_{12}$ and $G_{22}$ and their derivates are similar; we omit the details.

\noindent{\bf Case (iii).} Suppose now $\alpha_2\neq 0$. The idea of computing
$G_{11}$ is still the same as that in Case (ii). But now $c_0^+$ and $c_0^-$ in
(\ref{eq:exppx1}) and (\ref{eq:expqx1}) are zero, and
\[
  p_0^\pm = \left\{
    \begin{array}{ll}
\sqrt{\frac{\beta_1\pm \sqrt{\alpha_2}}{2c_{33}c_{22}}}, & \min_{\xi\in\mathbb{R}}\Delta(\xi^2)\notin (-\infty,0),\\
\sqrt{\frac{\beta_1\mp \sqrt{\alpha_2}}{2c_{33}c_{22}}}, & \min_{\xi\in\mathbb{R}}\Delta(\xi^2)\in (-\infty,0).\\
    \end{array}
  \right.
\]
Thus, it is easy to derive that, as $x\to 0$,
\begin{equation}
  \label{eq:G11:+-}
  G_{11}^{\pm}(x) = c_1^\pm(0)\log((p_0^\pm)^2x_2^2-x_1^2)+ {\cal O}\left(\sqrt{(p_0^\pm)^2x_2^2-x_1^2}\log((p_0^\pm)^2x_2^2-x_1^2)\right)
\end{equation}
exhibiting the logarithmic singularity of $G_{11}$. The analysis of $G_{12}$ and
$G_{22}$ is analogous to the previous case.

%
\section{Exact TBC for any exterior scattering problem}
In this section, we propose an exact TBC to truncate Navier's problem
(\ref{eq:gov}) in $\mathbb{R}^2\backslash\overline{D}$. Let $B(O,r)$ be of
sufficiently large radius $r>0$ to enclose $D$. Let $\Omega =
B(0,R)\backslash\overline{B(0,r)}$ for some constant $R>r$ and
$\Gamma=\partial\Omega$ be its boundary. Using Green's tensor $G$, we obtain
from the third Green's identity (c.f. Theorem 6.10 in \cite{mcl00}) that, for
any $x\in\Omega$,
\[
  u(x) = \int_{\Gamma}\left[[{\cal B}_{\nu(y)} G^{T}(y;x)]^{T} u(y) - G(y;x) {\cal
    B}_{\nu} u(y) \right]ds(y),\\
\]
where we recall that $\nu=[\nu_1,\nu_2]^{T}$ is the outer unit normal of the
boundary $\Gamma$ and ${\cal B}_\nu$ denotes the co-normal derivative operator
given by
\begin{align}
  {\cal B}_\nu =& [{\cal B}^{ij}_\nu]_{2\times 2}:=\nu_1(A_{11}\pd_1 +
                  A_{12}\pd_2) + \nu_2 (A_{21}\pd_1  + A_{22}\pd_2 ).
\end{align}
Inspired by SURC (\ref{eq:surc}), we impose a new SURC for the scattered elastic
wave $u$:
\begin{equation}
  \label{eq:rc}
  \lim_{R\to\infty} \int_{\partial B(0,R)}\left[[{\cal B}_{\nu(y)} G^{T}(y;x)]^{T} u(y) - G(y;x) {\cal
    B}_{\nu} u(y) \right]ds(y) = 0,
\end{equation}
uniformly for $x$ in any bounded domain. Its rigorous justification for $\rho\omega^2>0$ shall be presented
in a future work elsewhere. For dissipative media with ${\rm Im}(\rho\omega^2)>0$, (\ref{eq:rc}) is
evident as the dissipative Green tensor $G$ decays exponentially (c.f. section
3.3) and we must
look for an exponentially decaying solution $u$.

This radiation condition directly leads to Green's exterior
representation formula: for any $x\in \mathbb{R}^2\backslash\overline{B(0,r)}$,
\begin{equation}
  \label{eq:grep}
  u(x) = \int_{\partial B(0,r)}\left[[{\cal B}_{\nu(y)} G^{T}(y;x)]^{T} u(y) - G(y;x) {\cal
    B}_{\nu} u(y) \right]ds(y).
\end{equation}

By the standard jump relations (c.f. (7.5) in \cite{mcl00}),
we obtain an exact TBC
\begin{equation}
  \label{eq:tbc}
  \left( {\cal I}/2 - {\cal K}  \right) [u]= - {\cal S}[{\cal B}_\nu u],\quad{\rm on}\quad \partial B(0,r),
\end{equation}
where ${\cal I}$ denotes the identity operator, the single- and double-layer
potential operators ${\cal S}$ and ${\cal K}$ are defined by
\begin{align}
  {\cal S}[\phi] =& \int_{\partial B(0,r)} G(y;x)\phi(y) ds(y),\\
  {\cal K}[\phi] =& p.v.\int_{\partial B(0,r)} [{\cal B}_{ \nu(y) } G^{T}(y;x)]^{T}\phi(y) ds(y),
\end{align}
$\phi=[\phi_1,\phi_2]^{T}$, and $p.v.$ indicates Cauchy's principal value.

The TBC condition exactly truncates the unbounded scattering problem
(\ref{eq:gov}) onto the bounded domain $B(0,r)\backslash \overline{D}$ so that
any standard solvers can apply readily. In this paper, we shall derive a BIE
solver only for Navier's problem (\ref{eq:gov}) based on the new TBC
(\ref{eq:tbc}). To make the TBC applicable in practice, we require high-accuracy
quadrature rules to discretize the integrals in (\ref{eq:tbc}).
\subsection{High-accuracy discretization scheme}
In this subsection, we propose a high-accuracy quadrature rule to
discretize the two operators ${\cal S}$ and ${\cal K}$.
Consider the single-layer operator ${\cal S}$ first. Let $\partial B(0,r)$ be
parameterized by 
\[
  x(t) = (x_1(t),x_2(t))=(r\cos t,r\sin t), \quad 0\leq t\leq 2\pi. 
\]
Thus,
\begin{equation}
  \label{eq:dis:S}
  {\cal S}[\phi](x(t)) = r\int_{0}^{2\pi} G(x(t)-x(\tilde{t}))\phi(x(\tilde{t})) d\tilde{t}.
\end{equation}
Based on the discussion in section 4, it can be seen that the kernel tensor
\[
  G(x(t)-x(\tilde{t})) = {\cal O}(\log|t-\tilde{t}|), \quad {\rm as}\quad t-\tilde{t}\to 0.
\]
Thus, we can simply use Alpert's six-order quadrature rule
\cite{alp99,luluqia18} to discretize the integral in (\ref{eq:dis:S}) based on a
uniform discretization of $t$ over $[0,2\pi]$, i.e.
$\{t_j=\frac{2j\pi}{N}\}_{j=0}^{N-1}$, to approximate
\[
  {\cal S}[\phi]\left[
    \begin{matrix}
      x(t_0)\\
      x(t_1)\\
      \vdots\\
      x(t_{N-1})
    \end{matrix}
  \right]\approx
  \left[
    \begin{matrix}
      {\bm S}_{11} & {\bm S}_{12}\\
      {\bm S}_{21} & {\bm S}_{22}\\
    \end{matrix}
  \right]
  \left[
    \begin{matrix}
   {\bm \phi}_1\\   
   {\bm \phi}_2   
    \end{matrix}
  \right]
,
\]
where ${\bm S}_{ij}$ denotes an $N\times N$ matrix and ${\bm
  \phi}_j=[\phi_j(x(t_0)),\phi_j(x(t_1)),\cdots, \phi_j(x(t_{N-1}))]^{T}$.

We now turn to the discretization of the double-layer operator ${\cal K}$ given by
\begin{align}
  \label{eq:BvG}
  {\cal K}[\phi](x)
  =&p.v.\int_{\partial B(0,r)} \left[
    \begin{matrix}
      {\cal B}_{\nu(y)}^{11}G_{11} + {\cal B}_{\nu(y)}^{12}G_{21} & {\cal B}_{\nu(y)}^{21}G_{11} + {\cal B}_{\nu(y)}^{22}G_{21} \\
      {\cal B}_{\nu(y)}^{11}G_{12} + {\cal B}_{\nu(y)}^{12}G_{22}& {\cal B}_{\nu(y)}^{21}G_{12} + {\cal B}_{\nu(y)}^{22}G_{22}\\
    \end{matrix}
  \right](x;y) \phi(y)ds(y). 
\end{align}
To simplify the presentation, we consider the case $\alpha_2\neq 0$ only. 
Take
\[
 J(x):=\int_{\partial B(0,r)} {\cal B}_{\nu(y)}^{11}G_{11}(x;y)\phi_1 ds(y)
\]
as an example. Based on the splitting of $G$, considering the logarithmic
behavior of $G_{11}^\pm$ in (\ref{eq:G11:+-}), we rewrite
\begin{align*}
  {\cal B}_{\nu(y)}^{11} 
  =& f^\pm_1(y)[-(p_0^\pm)^2\nu_1(y)\partial_1 + \nu_2(y)\partial_2] + f_2^\pm(y)\partial_\tau,
\end{align*}
where $\tau$ denotes the unit tangential vector, and
\begin{align*}
  f^\pm_1(y) =\frac{-c_{11}\nu_1^2(y)-c_{33}\nu_2^2(y)}{(p_0^\pm)^2 \nu_1(y)^2 - \nu_2(y)^2},\quad
  f^\pm_2(y) =\frac{[(p_0^+)^2c_{33}+c_{11}]\nu_1(y)\nu_2(y)}{(p_0^\pm)^2 \nu_1(y)^2 - \nu_2(y)^2}.
\end{align*}
It can be verified that under $(C_0)$, $(p_0^\pm)^2$ are either
negative or complex so that $(p_0^\pm)^2 \nu_1(y)^2 -
\nu_2(y)^2\neq 0$ for any $y\in \partial D$. Therefore, $f_j^\pm(x(t))$ are analytic for
$t\in[0,2\pi]$. Thus,
\begin{align*}
  {\cal B}_{\nu(y)}^{11} G^\pm_{11}(x;y) =& c_0^{\pm}(0)f_1(y)[-(p_0^\pm)^2\nu_1(y)\partial_1 + \nu_2(y)\partial_2]\log((p_0^\pm)^2(x_2-y_2)^2-(x_1-y_1)^2)  \\
                                      &+ {\cal B}_{\nu(y)}^{11}[G^\pm_{11}(x;y) - c_0^{\pm}(0)\log((p_0^\pm)^2(x_2-y_2)^2-(x_1-y_1)^2)]\\
  &+ c_0^{\pm}(0)f^\pm_2(y)\partial_\tau \log((p_0^\pm)^2(x_2-y_2)^2-(x_1-y_1)^2).
\end{align*}
As $x\to y$ along $\partial B(0,r)$, the sum of the first two terms is
logarithmically singular by (\ref{eq:G11:+-}) and
\begin{align*}
  &[-(p_0^\pm)^2\nu_1(y)\partial_1 + \nu_2(y)\partial_2]\log((p_0^\pm)^2(x_2-y_2)^2-(x_1-y_1)^2)\\
   =&\frac{2(p_0^\pm)^2\nu_1(y)(y_1-x_1) + 2\nu_2(y)(p_0^\pm)^2(y_2-x_2)}{(p_0^\pm)^2(x_2-y_2)^2-(x_1-y_1)^2} = {\cal O}(1).
\end{align*}
On the other hand, 
\begin{align*}
  &\int_{\partial B(0,r)} f_{2}(y)\partial_\tau \log((p_0^+)^2(x_2-y_2)^2-(x_1-y_1)^2) \phi_{1}(y) ds(y)\\
  =& -  \int_{\partial B(0,r)} \log((p_0^+)^2(x_2-y_2)^2-(x_1-y_1)^2) \frac{d}{ds}[f_{2}(y)\phi_{1}(y)] ds(y).
\end{align*}
Thus, the integrand of $J(x)$ is transformed into the sum of integrals of
logarithmically singular integrands so that Alpert's quadrature applies
to approximate  
\[
  {\bm J}\approx
  {\bm A} {\bm \phi}_1 + {\bm B}^+{\bm \psi}^+_1 + {\bm B}^-{\bm \psi}^-_1,
\]
for three $N\times N$ matrices ${\bm A}$ and ${\bm B}^\pm$, where ${\bm J} =
[J(x(t_0)), J(x(t_1)), \cdots, J(x(t_{N-1}))]^{T}$, and
\begin{align*}
  {\bm \psi}^\pm_1=&\left[
    \begin{matrix}
      \frac{d}{dt}\left[\phi_1(x(t))f_2^\pm(x(t)) \right]|_{t=t_0}\\
      \frac{d}{dt}\left[\phi_1(x(t))f_2^\pm(x(t)) \right]|_{t=t_1}\\
      \vdots\\
      \frac{d}{dt}\left[\phi_1(x(t))f_2^\pm(x(t)) \right]|_{t=t_{N-1}}\\
    \end{matrix}
  \right] = {\bm T}^\pm_1 {\bm \phi}_1 + {\bm T}^\pm_2 {\bm \phi}'_1,\\
  {\bm T}^\pm_1 =& {\rm Diag}\{[f^\pm_2]'(x(t_0)),\cdots, [f^\pm_2]'(x(t_{N-1}))\},\\
  {\bm T}^\pm_2 =& {\rm Diag}\{f^\pm_2(x(t_0)),\cdots, f^\pm_2(x(t_{N-1}))\},\\
  {\bm \phi}'_1 =& {\rm Diag}\{\phi_1'(x(t_0)),\cdots, \phi_1'(x(t_{N-1}))\}.
\end{align*}
Moreover, considering that $\phi(x(t))$ is a $2\pi$-periodic smooth function for
$t\in[0,2\pi]$, we find, by the method of fast Fourier transform \cite{tre00}, a
spectral-accuracy $N\times N$ differentiation matrix ${\bm D}$ such that $ {\bm
  \phi}_1' = {\bm D}{\bm \phi}_1$. Consequently,
\[
  {\bm J} \approx [{\bm A} + {\bm B}^+({\bm T}^+_1+{\bm T}^+_2 {\bm D}) + {\bm
    B}^-({\bm T}^-_1+{\bm T}^-_2 {\bm D})]{\bm \phi}_1.
\]
One similarly discretizes the other integrals in (\ref{eq:BvG}) to approximate
\[
  {\cal K}[\phi]\left[
    \begin{matrix}
      x(t_0)\\
      x(t_1)\\
      \vdots\\
      x(t_{N-1})
    \end{matrix}
  \right]\approx
  \left[
    \begin{matrix}
      {\bm K}_{11} & {\bm K}_{12}\\
      {\bm K}_{21} & {\bm K}_{22}\\
    \end{matrix}
  \right]
  \left[
    \begin{matrix}
   {\bm \phi}_1\\   
   {\bm \phi}_2   
    \end{matrix}
  \right]
.
\]
Then, the TBC (\ref{eq:tbc}) approximately becomes
\begin{align}
  \label{eq:app:tbc}
\left[
    \begin{matrix}
      {\bm I}/2-{\bm K}_{11} & -{\bm K}_{12}\\
      -{\bm K}_{21} & {\bm I}/2-{\bm K}_{22}\\
    \end{matrix}
  \right]
  \left[
    \begin{matrix}
   {\bm u}_1\\   
   {\bm u}_2   
    \end{matrix}
  \right] \approx -\left[
    \begin{matrix}
      {\bm S}_{11} & {\bm S}_{12}\\
      {\bm S}_{21} & {\bm S}_{22}\\
    \end{matrix}
  \right]
  \left[
    \begin{matrix}
   {\bm v}_1\\   
   {\bm v}_2   
    \end{matrix}
  \right],
\end{align}
where ${\bm u}_j = [u_j(x(t_0)),\cdots, u_j(x(t_{N-1}))]^{T}$ and ${\bm v}_j =
[{\cal B}_{\nu}u_j(x(t_0)),\cdots, {\cal B}_{\nu}u_j(x(t_{N-1}))]^{T}$.

The above discretization scheme for the smooth circle $\partial B(0,r)$ can be
trivially extended to any smooth closed curves.
\subsection{Numerical evidences}
We now carry out numerical experiments to illustrate the correctness of the TBC
(\ref{eq:tbc}). All numerical algorithms in this paper are implemented in MATLAB
2020a on a 2019 Macbook Pro.

\begin{figure}[!ht]
  \centering
\includegraphics[width=0.45\textwidth]{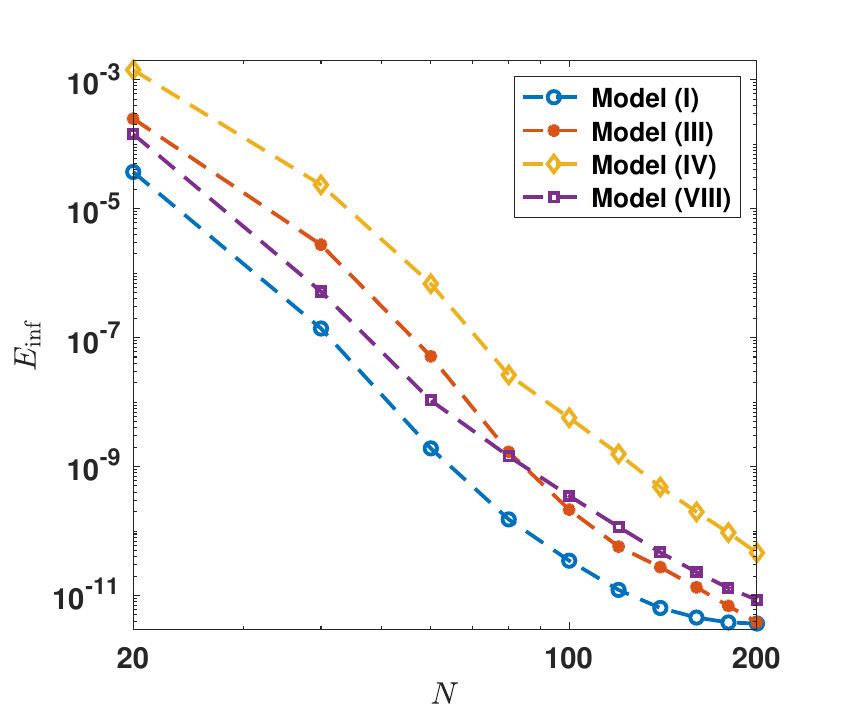}
\caption{Error $E_{\rm inf}$ for the approximate TBC (\ref{eq:app:tbc}) against $N$ for media (I), (III), (V) and (VIII).}
  \label{fig:err:N}
\end{figure}

Choose $r=1$, $N_0=8$, $M=20$, and let $u=G_1(x;x^*)$ be an elastic field
excited by the source point $x^*=(0.3,0)\in B(0,r)$. For all previous media in
section 3.3 with $\rho=\omega=1$, we choose proper values of $N$ such that the
following error quantity
\[
E_{\inf} = \left|\left|\left[
    \begin{matrix}
      {\bm I}/2-{\bm K}_{11} & -{\bm K}_{12}\\
      -{\bm K}_{21} & {\bm I}/2-{\bm K}_{22}\\
    \end{matrix}
  \right]
  \left[
    \begin{matrix}
   {\bm u}_1\\   
   {\bm u}_2   
    \end{matrix}
  \right] +\left[
    \begin{matrix}
      {\bm S}_{11} & {\bm S}_{12}\\
      {\bm S}_{21} & {\bm S}_{22}\\
    \end{matrix}
  \right]
  \left[
    \begin{matrix}
   {\bm v}_1\\   
   {\bm v}_2   
    \end{matrix}
  \right]\right|\right|_{\infty},
\]
is less than $10^{-9}$. The results are listed in Table~\ref{tab:tbc}.
\begin{table}
  \centering
  \begin{tabular}{c|c|c|c|c}
    \hline
   Medium & $N$ & $T_{\rm prep}$ & $T_{\rm cal}$ & $E_{\rm inf}$ \\\hline
    (I) & 80 & $14.2$s & $27.1$s & 1.5E$-10$\\\hline
    (II) & 240 & $10.3$s & $228.9$s & 5.8E$-10$\\\hline
    (III) & 100 & $12.3$s & $75.3$s & 2.1E$-10$\\\hline
    (IV) & 200 & $11.7$s & $219.5$s & 3.5E$-10$\\\hline
    (V) & 140 & $5.7$s & $134.9$s & 4.8E$-10$\\\hline
    (VI) & 140 & $11.3$s & $101.7$s & 2.3E$-10$\\\hline
    (VII) & 140 & $13.0$s & $97.0$s & 1.2E$-10$\\\hline
    (VIII) & 100 & $11.7$s & $50.3$s & 3.5E$-10$\\\hline
    (IX) & 100 & $14.5$s & $48.2$s & 3.4E$-10$\\\hline
  \end{tabular}
  \label{tab:tbc}
  \caption{Accuracy of the proposed TBC for all media in section 3.3.}
\end{table}
Here, $T_{\rm prep}$ records the running time for symbolically computing the
coefficients $c_n^\pm(x_2)$ in (\ref{eq:exppx1}) and (\ref{eq:expqx1}) for
$G_{11}$, and similar coefficients for the other components of $G$, and $T_{\rm
  cal}$ records the running time for evaluating $E_{\rm inf}$. Results in
Table~\ref{tab:tbc} provide a solid evidence for the correctness of the TBC,
regardless of the propagation behavior of wave fields at infinity. One may
increase $N$ to get more accurate approximations of the TBC. For example, we
show the relation of $E_{\rm inf}$ against $N$ for media (I), (III), (V) and
(VIII) in Figure~\ref{fig:err:N}, where $E_{\rm inf}$ decays at a rate roughly
proportional to $N^{-6}$, due to the order of accuracy of Alpert's rule we 
use.

\begin{figure}[!ht]
  \centering
(II)\includegraphics[width=0.215\textwidth]{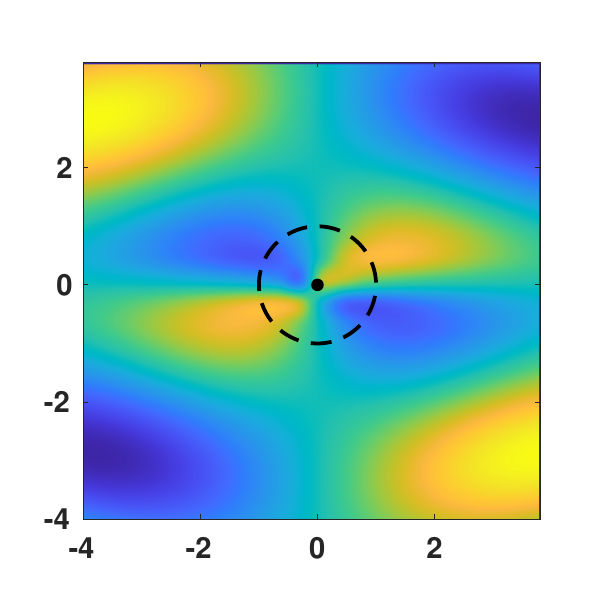}
\includegraphics[width=0.215\textwidth]{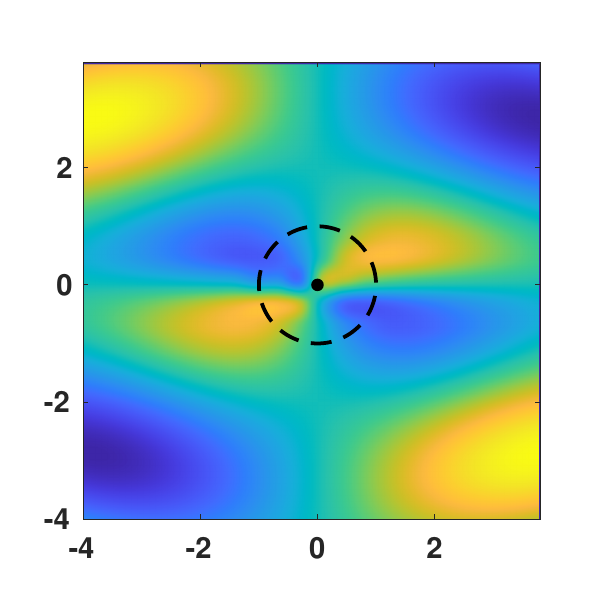}
(IV)\includegraphics[width=0.215\textwidth]{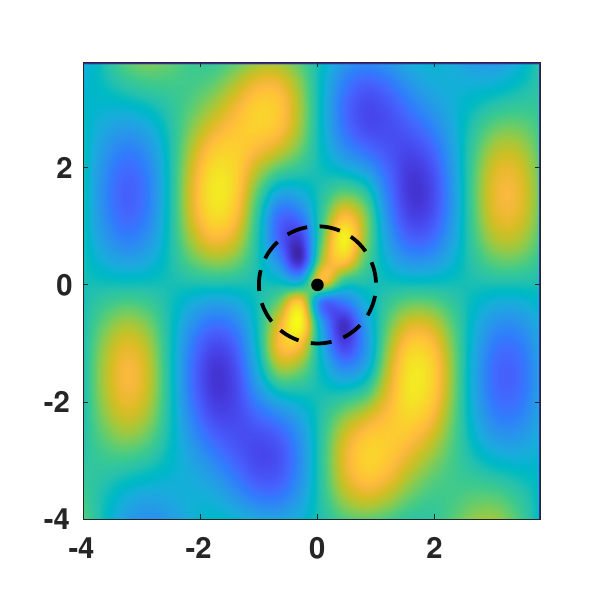}
\includegraphics[width=0.215\textwidth]{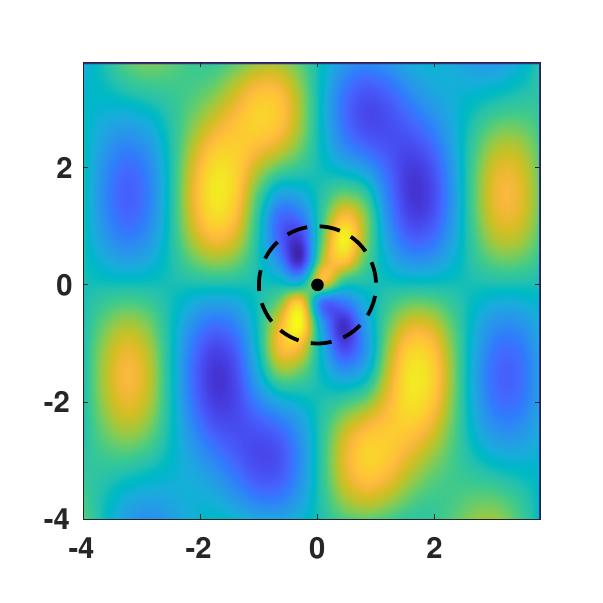}\\
(VII)\includegraphics[width=0.215\textwidth]{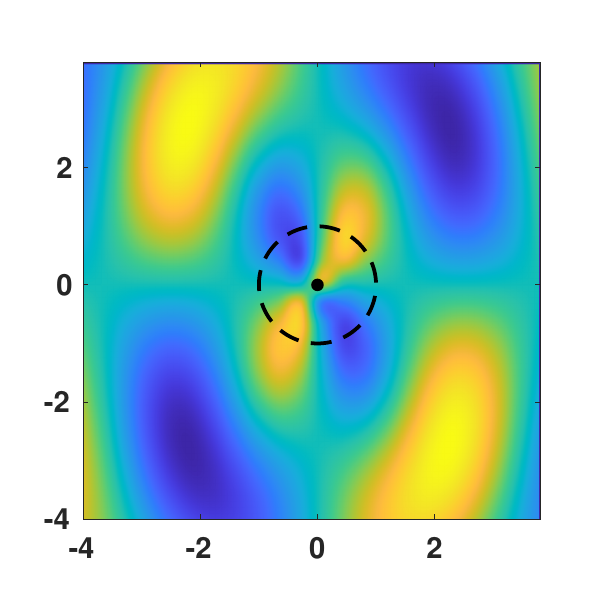}
\includegraphics[width=0.215\textwidth]{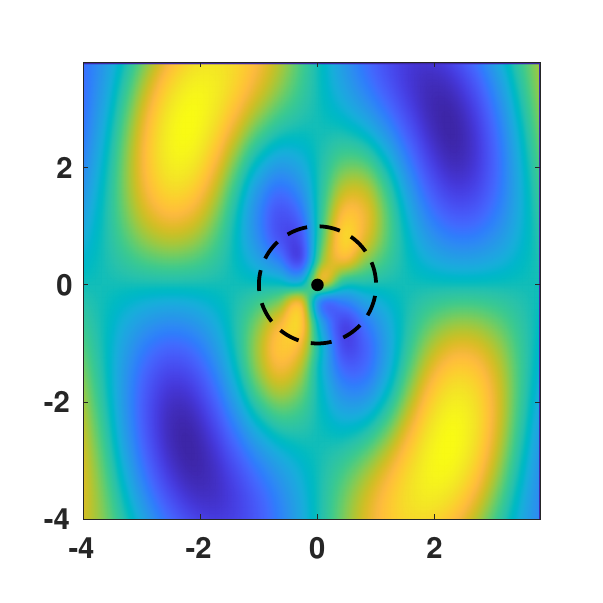}
(IX)\includegraphics[width=0.215\textwidth]{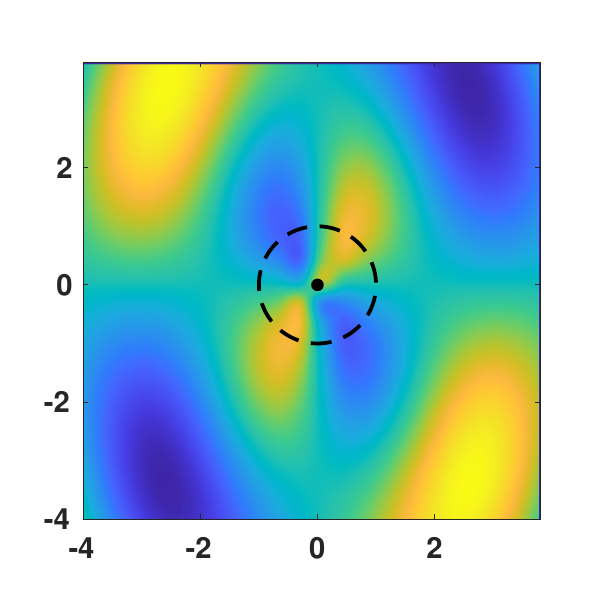}
\includegraphics[width=0.215\textwidth]{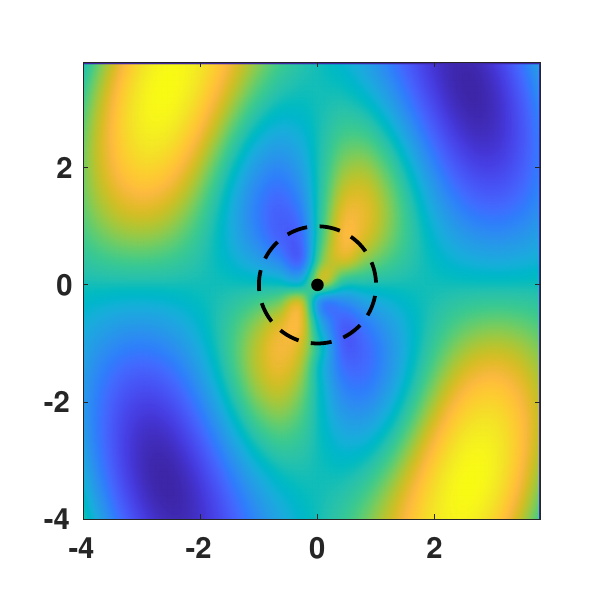}
\caption{Comparisons of numerical solutions (right) and exact solutions (left) for $G_{11}(x)$
  for media (II), (IV), (VII) and (IX) outside the TBC boundary $\partial
  B(0,1)$ indicated by the dashed lines. }
  \label{fig:grerep}
\end{figure}

If $u$ and ${\cal B}_{\nu} u$ on $\partial B(0,r)$ are available, then, in the
exterior domain $\mathbb{R}^2\backslash\overline{B(0,r)}$, $u$ can be computed
via Green's representation formula (\ref{eq:grep}). This is illustrated by
Figure~\ref{fig:grerep} for media (II), (IV), (VII) and (IX),
where numerical and exact solutions for $G_{11}(x)$ are compared, and where
dashed lines indicate the TBC boundary $\partial B(0,1)$. Numerical solutions
outside $B(0,1)$ are obtained by using exact values of $u$ and ${\cal B}_\nu u$
on $\partial B(0,1)$. The indistinguishability of numerical and exact solutions
provides another solid evidence for the correctness of the TBC (\ref{eq:tbc}).

\section{Typical scattering problems}

In this section, we consider two typical scattering problems to illustrate
the effectiveness of the proposed TBC. 

\subsection{Navier's Dirichlet-type problems}
The Dirichlet problem in section 1.3 can now be completely posed as
\begin{align}
  -\pd_1(A_{11}\pd_1 u + A_{12}\pd_2 u) -\pd_2(A_{21}\pd_1 u + A_{22}\pd_2 u) -\rho\omega^2 u =& 0,\quad {\rm in}\quad \mathbb{R}^2\backslash\overline{D},\\
  \label{eq:p1bc}
  u|_{\partial D} =& -u^{\rm inc}|_{\partial D},\\
  u\ {\rm satisfies}\ &\eqref{eq:rc}.
\end{align}
Here, the incident wave $u^{\rm inc}(x)$ can be a point-source wave
$G_j(x;x^*), j=1,2$, excited by a specified source point $x^*\not\in
\overline{D}$ or a plane wave $u^0 e^{\bi \sqrt{\rho}\omega x\cdot d}$, with a
nonzero polarization vector $u^0=[u^0_1,u^0_2]^{T}$ and a directional vector
$d=[d_1,d_2]^{T}$, determined in advance by
\begin{equation}
  \label{eq:plane:inc}
  \left[
    \begin{matrix}
      c_{11} d_1^2 + c_{33} d_2^2 &  (c_{12}+c_{33})d_1d_2\\
      (c_{12}+c_{33})d_1d_2 & c_{33}d_1^2 +c_{22}d_2^2 
    \end{matrix}
  \right]u^0 = u^0.
\end{equation}
To solve this unbounded domain problem,  we
directly make use of the TBC (\ref{eq:tbc}) on $\partial D$. Then, 
\begin{equation}
  \label{eq:method1}
  {\cal B}_{\nu} u = {\cal S}^{-1}({\cal I}/2-{\cal K}) [u^{\rm inc}],\quad {\rm
    on}\quad\partial D.
\end{equation}
Based on the high-accuracy discretization scheme in section 5.1, ${\cal B}_{\nu}
u$ can be accurately approximated. Once ${\cal B}_{\nu} u$ on $\partial D$
becomes available, Green's exterior representation formula (\ref{eq:grep})
applies to get $u$ in $\mathbb{R}^2\backslash\overline{D}$. One might argue for
an indirect approach to formulate the problem into a second-kind Fredholm
equation, but this is not the aim of the current paper.

\subsection{Navier's Neumann-type problems}
A Neumann problem replaces (\ref{eq:p1bc}) with 
\begin{align}
  \label{eq:p2bc}
  {\cal B}_{\nu} u|_{\partial D} =& -{\cal B}_{\nu}u^{\rm inc}|_{\partial D}.
\end{align}
It models wave scattered by an artificial elastic obstacle $D$ with the property
\[
  {\cal B}_{\nu}u^{\rm tot} = {\cal B}_{\nu}u + {\cal B}_{\nu}u^{\rm inc} =
  0,\quad{\rm on}\quad \partial D,
\]
where we recall $u^{\rm tot}$ denotes the total wave field $u+u^{\rm inc}$. Now,
the TBC (\ref{eq:tbc}) implies
\[
  u = ({\cal I}-{\cal K})^{-1}{\cal S}[{\cal B}_{\nu} u^{\rm inc}],
\]
and then Green's exterior representation formula (\ref{eq:grep}) applies to get
$u$ in $\mathbb{R}^2\backslash\overline{D}$.
\subsection{Numerical experiments}
We study three specific examples to illustrate the accuracy of the proposed approach.

Let $D$ be a kite-shaped domain with boundary $\partial D$
parameterized by 
\begin{equation}
  \label{eq:kite}
  x(t) = (\cos(2\pi t) + 0.65\cos(4\pi t) - 0.65, 1.5\sin(2\pi t)),\quad 0\leq
  t\leq 1.
\end{equation}
We consider two incident waves: (1) a plane wave $u_{\rm pw}^{\rm inc}(x) = A_0e^{\bi
  \sqrt{\rho}\omega x\cdot d}$ with a normalized vector $A_0$ and $d =
d_2[\sqrt{\frac{c_{33}}{c_{11}}} ,1]^{T}$ satisfying (\ref{eq:plane:inc}); (2) a
point-source wave $u_{\rm ps}^{\rm inc}(x) = G_1(x;x^*)$ with $x^* = (0,3)\notin D$.

In the implementation, we assume $\rho=\omega=1$, let $M=20$ in
(\ref{eq:G11:case2}) and then truncate the series after $n=8$ to evaluate
$G_{11}$ and, similarly, the other components in $G$ as well as their
derivatives. Define
\[
  {\rm Relative\ error} = \frac{||{\bm \phi}_{\rm num} - {\bm \phi}_{\rm ref}||_{\rm inf}}{||{\bm \phi}_{\rm ref}||_{\infty}},
\]
where ${\bm \phi}_{\rm num}$ represents a numerical solution of the co-normal
derivative ${\cal B}_{\nu} u$ for Dirichlet problems or that of the scattered
field $u$ for Neumann problems, at the $N$ grid points on $\partial D$, and
${\bm \phi}_{\rm ref}$ represents a reference solution, obtained for a
sufficiently large number $N$.

\begin{figure}[!ht]
  \centering
(a)\includegraphics[width=0.217\textwidth]{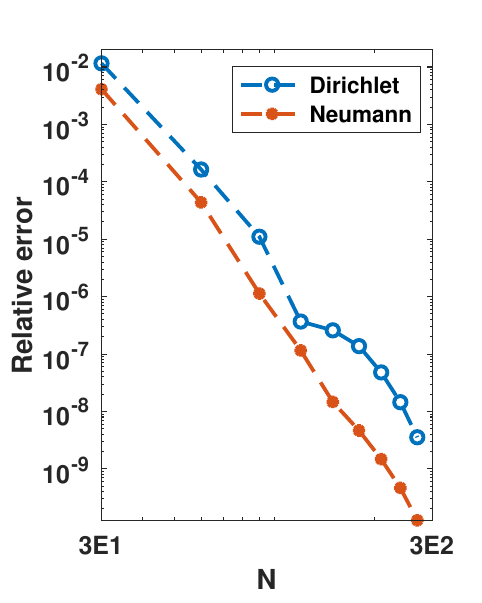}
(b)\includegraphics[width=0.267\textwidth]{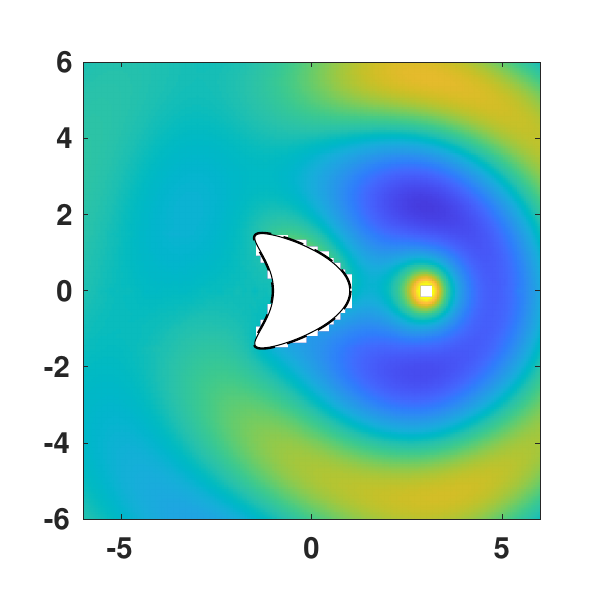}
(c)\includegraphics[width=0.267\textwidth]{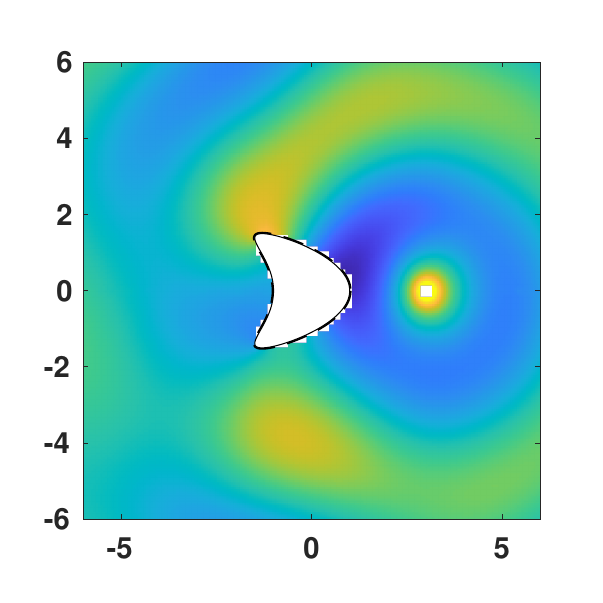}
\caption{Scattering problems for medium (III) with the point-source incidence
  $u^{\rm inc}_{\rm ps}$: (a) relative error vs. $N$; (b) real part of $u_1$ for the
  Dirichlet problem; (c) real part of $u_1$ for the Neumann problem.}
  \label{fig:M3:ps}
\end{figure}

First, we study the Dirichlet and Neumann problems for medium (III) by
specifying the point-source incident wave $u_{ps}^{\rm inc}$. Numerical results
are shown in Figure~\ref{fig:M3:ps},
where (a) shows relative error against $N$ with both
axes in logarithmic scale, and (b) and (c) show real parts of $u_1$ in
$[-6,6]^2$ for Dirichlet and Neumann problems, respectively.

\begin{figure}[!ht]
  \centering
(a)\includegraphics[width=0.217\textwidth]{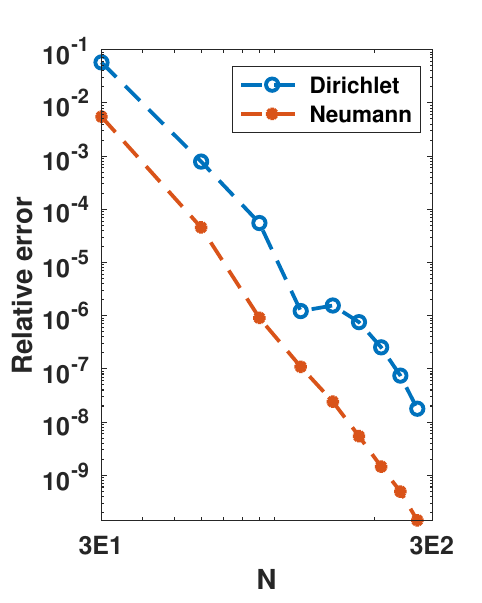}
(b)\includegraphics[width=0.267\textwidth]{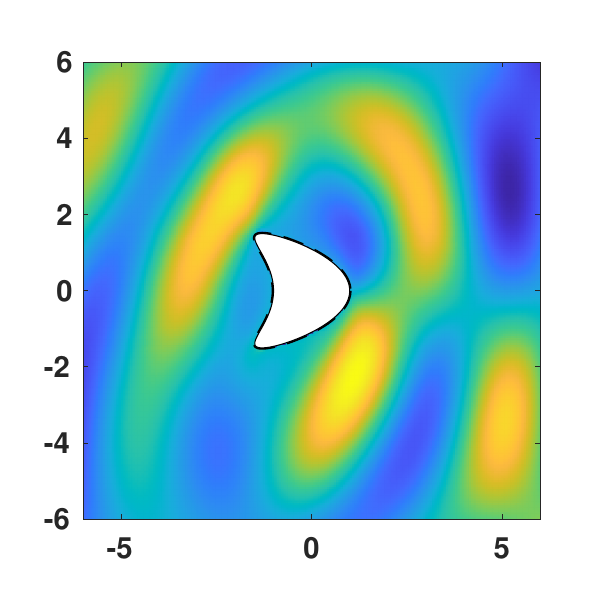}
(c)\includegraphics[width=0.267\textwidth]{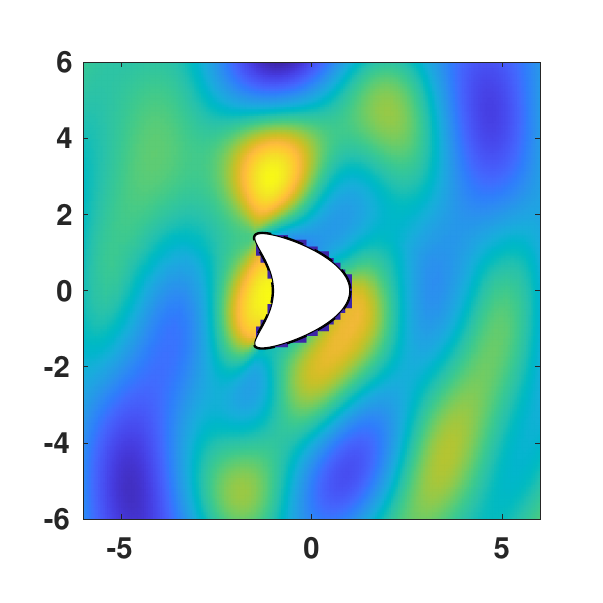}
\caption{Scattering problems for medium (VI) with the plane-wave incidence
  $u^{\rm inc}_{pw}$: (a) relative error vs. $N$; (b) real part of $u_1$ for the
  Dirichlet problem; (c) real part of $u_1$ for the Neumann problem.}
  \label{fig:M7:pw}
\end{figure}

Second, we study the Dirichlet and Neumann problems for medium (VI) by
specifying the
plane incident wave $u_{pw}^{\rm inc}$. Numerical results are shown in
Figure~\ref{fig:M7:pw},
where (a) shows relative error against $N$, and (b) and (c) show real parts of
$u_1$ in $[-6,6]^2$ for Dirichlet and Neumann problems, respectively.

\begin{figure}[!ht]
  \centering 
(a)\includegraphics[width=0.217\textwidth]{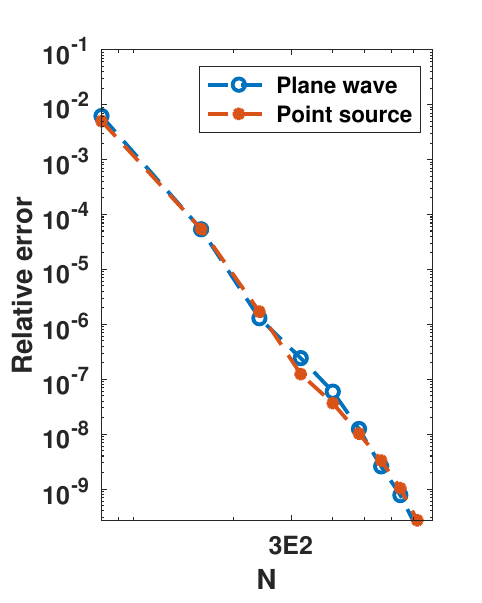}
 (b)\includegraphics[width=0.267\textwidth]{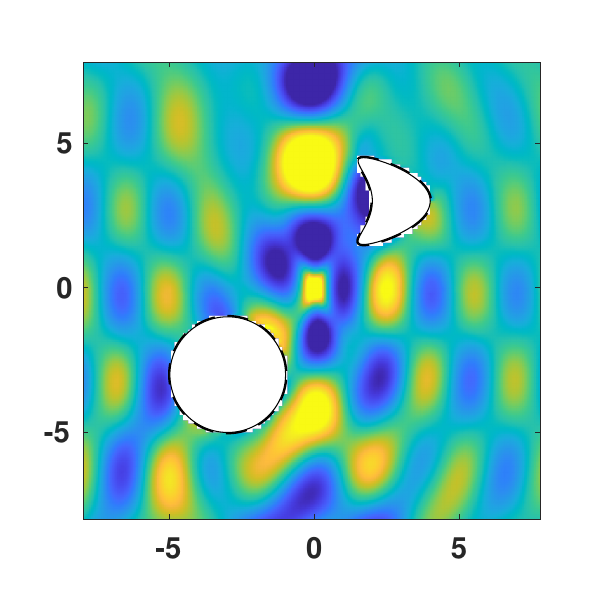}
 (c)\includegraphics[width=0.267\textwidth]{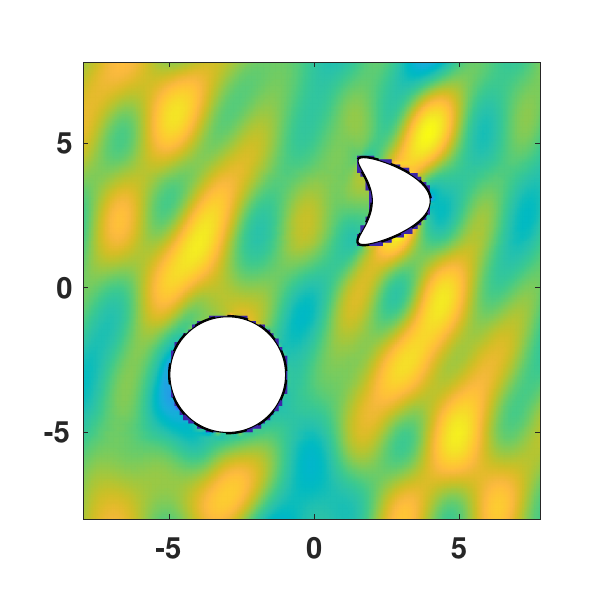}
 \caption{Neumann scattering problems for medium (IV): (a) relative error vs.
   $N$; (b) real part of $u_1$ for the point-source incidence $u^{\rm inc}_{\rm
     ps}$; (c) real part of $u_1$ for the plane-wave incidence $u^{\rm inc}_{\rm
     pw}$.}
  \label{fig:M4:neu}
\end{figure}

To conclude this section, in the last example, suppose $D$ consists of two well
separated domains: the previously introduced kite-shaped domain (\ref{eq:kite})
centered at $(3,3)$ and a disk of radius $2$ centered at $(-3,-3)$. Clearly, it
loses efficiency to use a big circle to enclose $D$, as it enlarges the
computational domain. Alternatively, we can directly establish the TBC
(\ref{eq:tbc}) on $\partial D$. We study the Neumann problem for medium (IV) for
both the point-source incidence $u^{\rm inc}_{\rm ps}$ and the plane-wave
incidence $u^{\rm inc}_{\rm pw}$. The results are shown in
Figure~\ref{fig:M4:neu}.
where (a) shows relative error against $N$, and (b) and (c) show real parts of
$u_1$ in $[-8,8]^2$ for the point-source incidence $u^{\rm inc}_{\rm ps}$ and
the plane-wave incidence $u^{\rm inc}_{\rm pw}$, respectively.

\section{Conclusions and Discussions}
We derived two-dimensional elastrodynamic Green's tensor in anisotropic media by
the method of Fourier transform and the use of SIP, and presented a rigorous
theory to completely classify the propagation behavior of Green's tensor for any
tensor $C$ satisfying $(C_0)$. Based on Green's tensor, we introduced a new SURC
(\ref{eq:rc}) to characterize elastic scattered waves at infinity due to bounded
scatterers. Based on (\ref{eq:rc}), the simple TBC (\ref{eq:tbc}) in terms of
single- and double-layer integral operators were used to analytically truncate
the scattered waves, whether backward waves exist or not. A number of numerical
experiments have validated the correctness of the TBC. Nevertheless, many
significant problems arise and deserve to be studied. We list a few but
certainly incomplete problems below.

Firstly, an essential problem is: how to rigorously justify the correctness of
the new SURC (\ref{eq:rc})? We believe that this can be solved by directly proving the well-posedness of the
scattering problem in $H^1_{\rm loc}(\mathbb{R}^2\backslash\overline{D})$ under
(\ref{eq:rc}). 

The second problem is: how to simplify the new SURC (\ref{eq:rc})? Based on the
relation of SRC (\ref{eq:src}) with SURC (\ref{eq:surc}), we believe that this
problem can be solved by carefully analyzing the asymptotic behavior of $G(x)$
for $|x|\gg 1$.

Considering that the exact TBC (\ref{eq:tbc}) is a nonlocal boundary condition, the third
problem is: whether a PML-like local boundary condition exists or not?
Discouragingly, the answer might be ``No'', as it is impossible to transform
waves such as $c_1e^{\bi k_1 x_2} + c_2e^{-\bi k_2 x_2}, k_j>0,$ to
exponentially decaying waves due to the completely different behavior at their
essential singularities in $\mathbb{C}$.

The fourth problem is: how to extend the current approach to more general media,
e.g., general anisotropic media with $c_{13}c_{23}\neq 0$, three-dimensional
general anisotropic media, etc.? Based on the same idea, the results must be as
simple as these stated in Theorems~\ref{thm:x2w}, \ref{thm:4.2.1}, and
\ref{thm:4.2.2} of the current paper, although the derivations are tolerably
tedious.

The last, most significant problem is: how to extend the current approach to
time-domain problems? An immediate approach might be Fourier inverse
transforming frequency-domain Green's tensor to get time-domain Green's tensor
so as to derive a time-domain TBC.

\bibliographystyle{plain}
\bibliography{wt}
\end{document}